\input amstex
\documentstyle{amsppt}
\magnification1200

\def\A{{\Cal A}}

\def\B{{\Cal B}}
\def\bx{{\boxed{\phantom{\square}}\kern-.4pt}} 

\def\C{{\Cal C}}
\def\CC{{\Bbb C}}
\def\Chi{{\operatorname{X}}}
\def\CMN{{\CC^N\!\ot\CC^M}}
\def\crc{{\raise.24ex\hbox{$\sssize\kern.1em\circ\kern.1em$}}}
\def\CSlm{{\CC\!\cdot\!S_{n+m}}}
\def\CSm{{\CC\!\cdot\!S_m}}
\def\CSn{{\CC\!\cdot\!S_n}}

\def\d{\partial}
\def\de{\delta}
\def\De{\Delta}
\def\dim{{\operatorname{dim\,}}}
\def\DMN{{\Cal{D}}}

\def\End{{\operatorname{End}}}
\def\EndCN{\End(\CC^N)}

\def\enddemos{{\ $\square$\enddemo}}

\def\Ga{\Gamma}
\def\Gaz{{\Ga(z_1\lc z_m)}}
\def\ge{\geqslant}
\def\glM{{\frak{gl}_M}}
\def\glMN{{\glN\times\glM}}
\def\glN{{\frak{gl}_N}}

\def\id{{\operatorname{id}}}
\def\io{\alpha}

\def\la{\lambda}
\def\lap{\lambda^{\ts\prime}}
\def\La{\Lambda}
\def\Lac{{\operatorname{\Lambda}^{\!c}}} \def\laps{{\lap_{s}}}
\def\lapsmin{{\lap_{s-1}}}
\def\Lar{{\operatorname{\Lambda}^{\!r}}} \def\le{\leqslant}

\def\mup{\mu^{\ts\prime}}

\def\Muc{{\operatorname{M}^c}}
\def\mups{{\mup_{s}}}
\def\mupsmin{{\mup_{s-1}}}

\def\om{\omega}
\def\ot{\otimes}

\def\PDMN{{\Cal{PD}}}
\def\PMN{{\Cal{P}}}
\def\phi{\varphi}
\def\phik{{\phi_{ijk}(u,v,w)}}
\def\Phiz{{\Phi_{\la}(z_1\lc z_n)}}

\def\R{{\Cal R}}
\def\ro{\rho}

\def\SSS{{\Cal S}}
\def\si{\sigma}
\def\slN{{\frak{sl}_N}}

\def\The{\Theta}
\def\Thez{{\The_\la(z_1\lc z_n)}}
\def\Tla{{\Cal T_\la}}

\def\tr{{\operatorname{tr}}}
\def\ts{{\hskip1pt}}

\def\UN{{\operatorname{U}(\glN)}}
\def\Up{\Upsilon}
\def\Upz{{\Up_\la(z_1\lc z_n)}}

\def\Xiu{{\operatorname{Z}(u)}}

\def\YN{{\operatorname{Y}(\glN)}}
\def\YslN{{\operatorname{Y}(\slN)}}

\def\Z{{\Cal Z}}
\def\ZN{{\operatorname{Z}(\glN)}}

\expandafter\ifx\csname bethe.def\endcsname\relax \else\endinput\fi
\expandafter\edef\csname bethe.def\endcsname{%
 \catcode`\noexpand\@=\the\catcode`\@\space}
\catcode`\@=11

\mathsurround 1.6pt

\def\hcor#1{\advance\hoffset by #1}
\def\vcor#1{\advance\voffset by #1}
\let\bls\baselineskip  \let\ignore\ignorespaces
\def\vsk#1>{\vskip#1\bls} \let\adv\advance 
\def\vv#1>{\vadjust{\vsk#1>}\ignore} \def\vvv#1>{\vadjust{\vskip#1}\ignore}
\def\vvn#1>{\vadjust{\nobreak\vsk#1>\nobreak}\ignore}
\def\vvvn#1>{\vadjust{\nobreak\vskip#1\nobreak}\ignore}
\def\setnormalbls{\edef\normalbls{\bls\the\bls}}
\def\setmaths{\edef\maths{\mathsurround\the\mathsurround}}

 \let\nt\noindent 
\def\nn#1>{\noalign{\vskip #1pt}} \def\NN#1>{\openup#1pt}
 
\let\Sum\sum \def\sum{\Sum\limits} 
\let\Prod\prod \def\prod{\Prod\limits} \let\Int\int \def\int{\Int\limits}

\let\=\m@th \def\&{.\kern.1em} \def\>{\!\;} \def\:{\!\!\;}

\ifx\plainfootnote\undefined \let\plainfootnote\footnote \fi
\expandafter\ifx\csname amsppt.sty\endcsname\relax
 
\else \fi

\newbox\s@ctb@x
\def\s@ct#1 #2\par{\removelastskip\vsk>
 \vtop{\bf\setbox\s@ctb@x\hbox{#1} \parindent\wd\s@ctb@x
 \ifdim\parindent>0pt\adv\parindent.5em\fi\item{#1}#2\strut}%
 \nointerlineskip\nobreak\vtop{\strut}\nobreak\vsk-.4>\nobreak}

\newbox\t@stb@x
\def\gadv{\global\advance} \def\gad#1{\gadv#1 1} 
\def\l@b@l#1#2{\def\n@@{\csname #2no\endcsname}%
\if *#1\gad\n@@ \expandafter\xdef\csname @#1@#2@\endcsname{\the\Sno.\the\n@@}%
 \else\expandafter\ifx\csname @#1@#2@\endcsname\relax\gad\n@@
 \expandafter\xdef\csname @#1@#2@\endcsname{\the\Sno.\the\n@@}\fi\fi}
\def\l@bel#1#2{\l@b@l{#1}{#2}\?#1@#2?}
\def\?#1?{\csname @#1@\endcsname}
\def\[#1]{\def\n@xt@{\ifx\t@st *\def\n@xt####1{{\setbox\t@stb@x\hbox{\?#1@F?}%
 \ifnum\wd\t@stb@x=0 {\bf???}\else\?#1@F?\fi}}\else
 \def\n@xt{{\setbox\t@stb@x\hbox{\?#1@L?}\ifnum\wd\t@stb@x=0 {\bf???}\else
 \?#1@L?\fi}}\fi\n@xt}\futurelet\t@st\n@xt@}
\def\(#1){{\rm\setbox\t@stb@x\hbox{\?#1@F?}\ifnum\wd\t@stb@x=0 ({\bf???})\else
 (\?#1@F?)\fi}}
\def\dff{\expandafter\d@f} \def\d@f{\expandafter\def}
\def\edff{\expandafter\ed@f} \def\ed@f{\expandafter\edef}

\newcount\Sno \newcount\Lno \newcount\Fno
\def\Section#1{\gadno\Fno=0\Lno=0\s@ct{\the\Sno.} {#1}\par} \let\Sect\Section
\def\section#1{\gad\Sno\Fno=0\Lno=0\s@ct{} {#1}\par} \let\sect\section
\def\l@F#1{\l@bel{#1}F} \def\<#1>{\l@b@l{#1}F} \def\l@L#1{\l@bel{#1}L}
\def\Tag#1{\tag\l@F{#1}} \def\Tagg#1{\tag"\llap{\rm(\l@F{#1})}"}
\def\Th#1{Theorem \l@L{#1}} \def\Lm#1{Lemma \l@L{#1}}
\def\Prop#1{Proposition \l@L{#1}}
\def\Cr#1{Corollary \l@L{#1}} \def\Cj#1{Conjecture \l@L{#1}}
 
\def\Proof#1.{\demo{\it Proof #1}}

\def\Par{\par\medskip} \def\setparindent{\edef\Parindent{\the\parindent}}
\def\Appendix{\Sno=64\let\p@r@\z@ %\parindent
 \def\Section##1{\gad\Sno\Fno=0\Lno=0 \s@ct{} \hskip\p@r@ Appendix
\\the\Sno
  \if *##1\relax\else {.\enspace##1}\fi\par} \let\Sect\Section
 \def\section##1{\gad\Sno\Fno=0\Lno=0 \s@ct{} \hskip\p@r@ Appendix%
  \if *##1\relax\else {.\enspace##1}\fi\par} \let\sect\section
 \def\l@b@l##1##2{\def\n@@{\csname ##2no\endcsname}%
 \if *##1\gad\n@@
 \expandafter\xdef\csname @##1@##2@\endcsname{\char\the\Sno.\the\n@@}%
 \else\expandafter\ifx\csname @##1@##2@\endcsname\relax\gad\n@@
 \expandafter\xdef\csname @##1@##2@\endcsname{\char\the\Sno.\the\n@@}\fi\fi}}

\let\le\leqslant \let\ge\geqslant
  \let\8\infty

\let\=\m@th  \def\_#1{_{\rlap{$\ssize#1$}}}

\def\lc{{,\ldots\hskip-0.1pt,\hskip1pt}}

\def\E(#1){\mathop{\hbox{\rm End}\,}(#1)} 
\def\id{\hbox{\rm id}}  
\def\tr{\hbox{\rm tr}}

\def\1{^{-1}} \def\vst#1{{\lower2.1pt\hbox{$\bigr|_{#1}$}}}

\let\logo@\relax
\let\m@k@h@@d\makeheadline \let\m@k@f@@t\makefootline
\def\makeheadline{\ifnum\pageno=1\headline={\hfil}\fi\m@k@h@@d}
\def\makefootline{\ifnum\pageno=1\footline={\hfil}\fi\m@k@f@@t}

\font\bigbf=cmbx10 scaled 1200
\centerline{\bigbf Yangians and Capelli Identities} \bigskip
\centerline{\smc Maxim Nazarov}
\bigskip\bigskip
\centerline{\it To Professor Alexandre Kirillov on the occasion
of his sixtieth birthday}
\bigskip\bigskip

\centerline{\S\ts\bf1.\ Introduction}\section{\,} \kern-20pt

\nt
In this article we apply representation theory of Yangians to the classical
invariant theory.
Let us consider the action of the Lie algebra $\glMN$ in the space $\Cal P$
of polynomial functions on $\CMN$.
This action is multiplicity-free, and its irredu\-cible components are
parametrized by Young diagrams $\la$ with not more than $N,M$ rows.
As a result, the space $\Cal{I}$
of $\glMN$\ts-\ts invariant differential operators on $\CMN$ with
polynomial coefficients splits into a direct sum of one-dimensional
subspaces parametrized by the diagrams $\la$.
It is easy to describe these subspaces.

Let $x_{ia}$ with $i=1\lc N$ and
$a=1\lc M$ be the standard coordinates on the vector space $\CMN$.
Let $\d_{ia}$ be the partial derivation with respect to the coordinate
$x_{ia}$. Suppose that the diagram $\la$ consists of
$n$ boxes. Let $\chi_\la$ be the irreducible character of the symmetric
group $S_n$ parametrized by $\la$.
Then the one-dimensional subspace in $\Cal I$ corresponding to $\la$
is spanned by the operator
$$
%c_{\ts\la}\,=
\sum_{\si\ts\in\,S_n}
\ \,
\sum_{i_1\lc i_n}
\,
\sum_{a_1\lc a_n}
\
\frac{\chi_\la\ts(\si)\!}{n!}
\cdot
x_{i_1a_1}\ldots x_{i_na_n}
\!\!\cdot
\d_{i_{\si(1)}a_1}\ldots\d_{i_{\si(n)}a_n} \Tag{1.1}
$$
where the indices $i_1\lc i_n$ and $a_1\lc a_n$ run through $1\lc N$
and $1\lc M.$

On the other hand, the action of the Lie algebra $\glN$ in the space $\Cal P$
extends to the action of the universal enveloping algebra $\UN$ by the
differential
operators with polynomial coefficients. The space $\Cal I$ 
is then the image of the centre of $\UN$;
see for instance [HU].
In this article we give an explicit formula for a central element of $\UN$ 
corresponding to the operator \(1.1). In the case when the diagram $\la$
has only one column this formula was discovered by A.\,Capelli [C].
In the other particular case when $\la$ consists of only one row 
this formula was found in [N1].

To derive an explicit formula for the general diagram $\la$ we
employ representation
theory of the Yangian $\YN$ of the Lie algebra $\glN$. The Yangian $\YN$
is a canonical deformation of the universal enveloping algebra
$\operatorname{U}(\glN[z])$ in the class of Hopf algebras [D1].
Moreover, 
%the algebra $\YN$
it contains $\UN$ as a subalgebra and admits a homomorphism
$\pi:\ts\YN\to\UN$ identical on $\UN$. Thus the irreducible representation
$\pi_\la$
of the Lie algebra $\glN$ corresponding to $\la$ can be regarded as
a representation of the algebra $\YN$. 

We use the notion of a universal $R$-matrix for the Hopf algebra $\!\YN,$
cf.~[D1]. Let $\Delta:\YN\to\YN\ot\YN$ be the comultiplication of $\YN$.
\text{Denote} by $\Delta^\prime$ the composition of $\Delta$
with permutation of tensor factors in $\YN\ot\YN$.
The algebra $\YN$ has a canonical family of automorphisms $\tau_z$ 
parametrized by $z\in\CC$.
The universal $R$-matrix for $\YN$
is a formal power series $\R(z)$ in $z\1$ with coefficients from $\YN\ot\YN$ 
and the leading term~$1$~such~that $$
\R(z)\cdot\id\ot\tau_z\bigl(\De^\prime(Y)\bigr) =
\id\ot\tau_z\bigl(\De(Y)\bigr)\cdot\R(z)\ts, \,\quad
Y\in\YN.
\Tag{1.0}
$$

For the description of the series $\R(z)$ see Section 3.
The image $\pi_\la\ot\pi\ts\bigl(\R(z)\bigr)$ is a rational function 
in $z$ with at most $n$ poles. These poles are contained in the collection
of the contents $c_1\lc c_n$ of the diagram $\la$; see Section 2. 
Denote by $\psi_\la$ the trace of the representation $\pi_\la$ of $\YN$ 
and consider the polynomial in~$z$ $$
\psi_\la\ot\pi\bigl(\R(z)\bigr)\cdot(z-c_1)\ldots(z-c_n) \Tag{1.2}
$$
valued in the algebra $\UN$; cf. [D2\ts,\ts RS]. 
We prove that the value of polynomial \(1.2) at $z=0$ 
is the central element corresponding to the operator \(1.1). 

There is an explicit formula for the polynomial \(1.2).
%Let $U_\la$ be the irreducible $S_n$-module corresponding to $\la$.
Denote by $\Lac$ the Young tableau obtained by filling the boxes of the
diagram $\la$ with the numbers $1\lc n$ by columns. 
Suppose that $c_1\lc c_n$ are the contents of the respective boxes of $\la$.
Denote by $S_\la$ and $T_\la$ the
subgroups in $S_n$ preserving the collections of numbers appearing 
respectively in every row and column of the tableau $\Lac$. 
Consider the element of the group ring $\CSn$ $$
\sum_{\si\in S_\la}
\,
\sum_{\ro,\ro^\prime\in T_\la}
\ \ro\ts\si\hskip-1pt\ro^{\ts\prime}
\cdot
\operatorname{sgn}(\ro\ro^{\ts\prime})
\cdot
y\ =
\sum_{\si\in S_n}
\ y_\si\hskip-1pt\cdot\si\,,
\ \quad
y_\si\in\CC
$$
where non-zero factor $y\in\CC$ is chosen to make this element of $\CSn$ 
an idempotent.

The Lie algebra $\glN$ acts in space $\CMN$ by linear combinations of the
operators $$
\sum_{1\le a\le M}\,x_{ia}\,\d_{ja}\,;
\qquad
i,j=1\lc N.
$$
The differential operator corresponding to the value of \(1.2) at $z=0$ 
is then
$$
%c_{\ts\la}\,=
\sum_{\si\ts\in\,S_n}
\ \,
\sum_{i_1\lc i_n}
\,
\sum_{a_1\lc a_n}
\
{y_\si}
\cdot\ts
\prod_k^\rightarrow\
\bigl(\ts x_{i_ka_k}\d_{i_{\si(k)}a_k}-c_k\cdot\de_{i_ki_{\si(k)}}\bigr) 
\Tag{1.3}
$$
where the index $k$ runs through $1\lc n$ and factors in the ordered product
are arranged from the left to right while $k$ increases. 
Here $\de_{ij}$ is the Kronecker delta. 

The equality of the differential operators \(1.1) and \(1.3)
has been proved in [O1].
It can be also verified by using the Jucys-Murphy elements in the group
ring of $S_n$; see forthcoming article [O2].
The present article contains a proof of that equality going back
to the origin \(1.2) of the formula \(1.3).
This proof is based on the estimation of the order of the pole at $z=0$
of~the rational function $\pi_\la\ot\pi_\mu\ts\bigl(\R(z)\bigr)$
for diagrams $\mu$ with less than $n$ boxes; cf.~[OO,\ts S].
That estimation is performed in Section 4. 
We employ the notion of a fusion procedure for the symmetric group $S_n$ 
introduced by I.\,Cherednik~[C1]. 
In Section 2 we give a concise account of the
fusion procedure; 
see [JKMO] for another exposition of relevant results from [C1].

I am very grateful to I.\,Cherednik for numerous illuminating conversations.
I am also grateful to F.\,Knop, M.\,Noumi and V.\,Tolstoy 
for stimulating discussions and valuable remarks. 
This work is a part of the project on representation theory of 
Yangians which has been started at the seminar of A.\,Kirillov in Moscow. 
I thank all participants of that seminar for providing a favourable 
atmosphere for the start.
I am especially indebted to A.\,Okounkov and G.\,Olshanski. 
Their results from [O1,OO] have inspired the present work. 

\newpage

\centerline{\S\ts\bf2.\ Fusion procedure for the symmetric group}\section{\,}
\kern-20pt

\nt
We will start with recalling several classical facts about irreducible
modules of the symmetric group $S_n$ over the complex field. 
They are parametrized by the Young diagrams $\la$ with exactly $n$ boxes.
We will denote by $U_\la$ the irreducible module of $S_n$ 
corresponding to the diagram $\la$. 
Consider the chain of subgroups
$$
S_1\subset S_2\subset\ldots\subset S_n 
$$ 
with respect to the standard embeddings. 
There is a canonical decomposition of the space $U_\la$ 
into the direct sum of
one-dimensional subspaces associated with this chain. 
These subspaces are parametrized by the {\it Young tableaux} of shape $\la$.
Each of these tableaux is a bijective filling of the boxes of $\la$ 
with numbers $1\lc n$ such that in every row and column the numbers 
increase from the left to the right and from the top to the bottom 
respectively. Denote by $\Cal T_\la$
the set of these tableaux.

For every tableau $\La\in\Tla$ denote by $U_\La$ the corresponding 
one-dimensional subspace in $U_\la$.
For any $m\in\{1\lc n\}$ consider the tableau obtained from $\La$ 
by removing each of numbers $m+1\lc n$. Let the Young diagram $\mu$ be 
its shape. Then $U_\La$ is contained in an irreducible $S_m$-submodule 
of $U_\la$ corresponding to $\mu$. Any basis of $U_\la$ formed by vectors 
$u_\La\in U_\La$ is called a {\it Young basis}. Let us fix an $S_n$-invariant 
inner product $\langle\,\ts\hskip0.5pt,\,\rangle$ on $U_\la$. The subspaces 
$U_\La$ are then pairwise orthogonal. We will be assuming that 
$\langle u_\La\ts,u_\La\rangle=1$ for each tableau $\La\in\Tla$.

Consider the {\it column tableau} of the shape $\la$ obtained by filling 
the boxes of $\la$ with $1\lc n$ by columns from the left to the right, 
downwards in each column. We will denote by $\Lac$ this tableau. 
Consider the diagonal matrix element of the $S_n$-module $U_\la$ 
corresponding to the vector $u_\Lac$
$$
\Phi_\la=\sum_{\si\in S_n}\ts
\langle\ts\si\ts u_\Lac\ts,u_\Lac\ts\rangle \!\cdot\hskip-1pt\si\1\,
\in\ \CSn\,.
\Tag{2.0}
$$
We will make use of the explicit formula for this matrix element contained 
in [Y1]. Denote by $S_\la$ and $T_\la$ the subgroups in $S_n$ preserving 
the collections of numbers appearing respectively in every row and column of 
the tableau $\Lac$. Put
$$
P_\la=\sum_{\si\in S_\la}\si\,,
\quad
Q_\la=\sum_{\si\in T_\la}\operatorname{sgn}(\si)\!\cdot\!\si\,. 
$$
As usual let $\lap_1\ts,\lap_2\ts,\ldots$ be the numbers of boxes in the 
columns of the diagram $\la$. Then 
$$
\Phi_\la=\frac{Q_\la\ts P_\la\ts Q_\la}
{\,\lap_1\ts!\,\lap_2\ts!\,\ldots\,}\ts. 
$$

We will also need for the matrix element $\Phi_\la$ an expression of a 
different kind.
For any distinct $i,j=1\lc n$ let $(i\ts j)$ be the transposition in 
the symmetric group $S_n$. Consider the rational function of two complex
variables $u,v$
valued in the group ring $\CSn$
$$
\phi_{ij}(u,v)=1-\frac{(i\ts j)}{u-v}.
$$
As a direct calculation shows, this rational function satisfies the 
equations
$$
\phi_{ij}(u,v)\ts\ts\phi_{ik}(u,w)\ts\ts\phi_{jk}(v,w)= 
\phi_{jk}(v,w)\ts\ts\phi_{ik}(u,w)\ts\ts\phi_{ij}(u,v) \Tag{2.1} 
$$
for all pairwise distinct $i,j,k$.
Evidently, for all pairwise distinct $i,j,k,l$ we have 
$$
\phi_{ij}(u,v)\ts\ts\phi_{kl}(z,w)=\phi_{kl}(z,w)\ts\ts\phi_{ij}(u,v) 
\Tag{2.2}
$$
Consider the rational function of $u,v,w$ appearing at either side of 
\(2.1). Denote by $\phik$ this function. The factor $\phi_{ik}(u,w)$ in 
\(2.1) has a pole at $u=w$. However, we have the following lemma.

\proclaim{Lemma 2.1}
Restriction of $\phik$ to the set of $(u,v,w)$ such that $u=v\pm1\ts,$ is 
regular at $u=w$. 
\endproclaim

\demo{Proof}
As a direct calculation shows, the restriction
$$
\phi_{ijk}(v\pm1,v,w)=
\Bigl(1\mp(i\ts j)\Bigr)\!\cdot\!\left(1-\frac{(i\ts k)+(j\ts k)}{v-w}\right).
$$
The latter rational function of $v,w$ is manifestly regular at $w=v\pm1$
\enddemos

\nt
For each $i\in\{1\lc n\}$ denote
$c_i=s-t$ if the number $i$ appears in the $s$-th column and the $t$-th row
of the tableau $\Lac$. The difference $s-t$ is then called the {\it content}
of the box of the diagram $\la$ occupied by the number $i$. 
For each $i$ let $z_i$ be a complex parameter. Equip the set of all pairs 
$(i,j)$ where $1\le i<j\le n$ with the lexicographical ordering. 
Introduce the ordered product over this set
$$
\prod_{(i,j)}^\rightarrow\ts\ts
\phi_{ij\ts}(c_i+z_i\ts,c_j+z_j).
\Tag{2.3}
$$
Consider this product as a rational function of the parameters $z_1\lc z_n$ 
valued in $\CSn$.
Denote by $\Phiz$ this rational function. Denote by $\Z$ the set of all 
tuples $(z_1\lc z_n)$ such that $z_i=z_j$ whenever the numbers $i$~and~$j$ 
appear in the same row of the tableau $\Lac$. The following
theorem goes back to [C1]~and~[J].

\proclaim{Theorem 2.2}
Restriction of $\Phiz$ to $\Z$ is regular at $z_1=\ldots=z_n$. 
The value of this restriction at $z_1=\ldots=z_n$ coincides with
\text{the matrix~element~$\hskip-1pt\Phi_\la\hskip-1pt$.} 
\endproclaim

\nt
We will present the main steps of the proof as separate propositions.
This proof follows [N2] and is based on Lemma 2.1. 
Another proof is contained in [JKMO].

\proclaim{Proposition 2.3}
Restriction of $\Phiz$ to $\Z$ is regular at $z_1=\ldots=z_n$.
\endproclaim

\demo{Proof}
We shall provide an expression for
the restriction of the function \(2.3) to $\Z$ which is manifectly regular
at $z_1=\ldots=z_n$.
Let us reorder the pairs $(i,j)$ in the product \(2.3) as follows. 
This reordering will not affect the value of the product due to the 
relations \(2.1) and \(2.2). 
Let $\C$ be the sequence of numbers obtained by reading the tableau 
$\Lac$ in the usual way, that is by rows from the top to the bottom, 
eastwards in every row. For~each $j\in\{1\lc n\}$ denote by $\A_j$ and 
$\B_j$ the subsequences of $\C$ consisting of all numbers $i<j$ which 
appear respectively before and after $j$ in that sequence. 
Now set $(i,j)\prec(k,l)$ if one of the following conditions is satisfied:

\medskip
\vbox{
\line{\ -\ the number $i$ appears in $\A_j$ while $k$ appears in $\B_l$;
\hfill}
\line{\ -\ the numbers $i$ and $k$ appear respectively in $\A_j$ and $\A_l$
where $j<l$;
\hfill}
\line{\ -\ the numbers $i$ and $k$ appear respectively in $\B_j$ and $\B_l$
where $j>l$;
\hfill}
\line{\ -\ we have the equality $j=l$ and $i$ appears before $k$ in $A_j$ or 
$B_j$.
\hfill}
}
\smallskip

{}From now on we assume that the factors in \(2.3) corresponding to the
pairs $(i,j)$ are arranged with respect to this new ordering. 
The factor $\phi_{\ts i\ts j\ts}(c_i+z_i\ts,c_j+z_j)$ has a pole at 
$z_i=z_j$ if and only if
the numbers $i$ and $j$ stand on the same diagonal of the tableau $\Lac$. 
We will then call the pair $(i,j)$ {\it singular.}
Note that the number $i$ occurs in the subsequence $\B_j$exactly when 
$i$ stands to the left and below of $j$ in the tableau $\Lac$.
In this case $c_j-c_i>1$ and the pair $(i,j)$ cannot be singular. 

Let a singular pair $(i,j)$ be fixed. Suppose that the number $i$ appears
in the $s$-th column and the $t$-th row of the tableau $\Lac$. 
In our new ordering the pair next after $(i,j)$ is $(h,j)$ 
where the number $h$ appears in the $(s+1)$-th column and the $t$-th row of
$\Lac$. In particular, we have $c_i=c_j=c_h-1$. Moreover, $(i,h)\prec(i,j)$.
Due to the relations \(2.1)\ts,\(2.2) the product 
$$
\prod_{(k,l)\prec(i,j)}^\rightarrow\ts
\phi_{kl\ts}(c_k+z_k\ts,c_l+z_l)
$$
is divisible on the right by
$\phi_{ih\ts}(c_i+z_i\ts,c_h+z_h)$. The restriction of the latter function to 
$z_i=z_h$ is just $1+(i\ts h)$. Note that the element 
$\bigl(1+(i\ts h)\bigr)/\ts2$
is an idempotent.

Now for each singular pair $(i,j)$ let us replace the two adjacent factors in 
\(2.3)
$$
\phi_{ij\ts}(c_i+z_i\ts,c_j+z_j)
\,
\phi_{hj\ts}(c_h+z_h\ts,c_j+z_j)
$$
by
$$
\gather
\phi_{ih\ts}(c_i+z_i\ts,c_h+z_h)
\,
\phi_{ij\ts}(c_i+z_i\ts,c_j+z_j)
\,
\phi_{hj\ts}(c_h+z_h\ts,c_j+z_j)/2
\\
=\phi_{ihj\ts}(c_i+z_i\ts,c_h+z_h,c_j+z_j)/\ts2. 
\endgather
$$
This replacement does not affect the value of the restriction to $\Z$ of the
function \(2.3). But the restriction to $z_i=z_h$ of the function 
$\phi_{ihj\ts}(c_i+z_i\ts,c_h+z_h,c_j+z_j)$ is regular at $z_i=z_j$ 
by Lemma 2.1
\enddemos

\nt
In our new ordering we have the decomposition $$
\Phiz=\Upz\cdot\Thez
$$
where $\Upz$ and $\Thez$ are products of the factors in \(2.3) which 
correspond
to the pairs $(i,j)$ with $i$ appearing in $\A_j$ and $\B_j$ respectively. 
The function $\Thez$ is regular at $z_1=\ldots=z_n$. Moreover, 
the value $\The$ of this function at $z_1=\ldots=z_n$ is 
invertible in $\CSn$. 
Denote by $\Phi$ and $\Up$ the values at $z_1=\ldots=z_n$ of the restrictions
to $\Z$ of $\Phiz$ and $\Upz$ respectively. 
Thus we have the equality $\Phi=\Up\ts\The$. 
We will prove that $\Phi=\Phi_\la$.

\proclaim{Proposition 2.4}
Let the numbers $k$ and $l$ stand next to each other in a row of the 
tableau $\Lac$.
Then the element $\Up$ in $\CSn$ is divisible on the right by $1+(k\ts l)$.
\endproclaim

\demo{Proof}
It suffices to assume that $k<l$. Then
due to the relations \(2.1) and \(2.2) the product $\Upz$ is divisible
on the right by $\phi_{kl\ts}(c_k+z_k\ts,c_l+z_l)$. 
The restriction of the latter function to $z_k=z_l$ is exactly the 
element $1+(k\ts l)$
\enddemos

\nt
Let $\io$ be the involutive antiautomorphism of the group ring $\CSn$ 
defined by setting $\io(g)=g\1$ for $g\in S_n$. 
By the relations \(2.1) and \(2.2) the product~\(2.3) is invariant 
with respect to this
antiautomorphism. So is the value $\Phi$ of its restriction to $\Z$. 

\proclaim{Corollary 2.5}
Let the numbers $k<l$ stand next to each other in the first row of 
the tableau $\Lac$.
Then the element $\Phi$ in $\CSn$ is divisible on the right by 
$$
\phi_{\ts kl\ts}(c_k\ts,c_l)\ts\ts\cdot\hskip-2pt 
\prod_{k<m<l}^\rightarrow\ts
\phi_{\ts ml\ts}(c_m\ts,c_l).
\Tag{2.5}
$$
The element $\Phi$ is then also divisible on the left by 
$$
\prod_{k<m<l}^\leftarrow\ts
\phi_{\ts ml\ts}(c_m\ts,c_l)
\cdot
\phi_{\ts kl\ts}(c_k\ts,c_l).
\Tag{2.6}
$$
\endproclaim

\demo{Proof}
By Proposition 2.4 the element $\Up$ is divisible on the right by 
$\phi_{\ts kl\ts}(c_k\ts,c_l).$ But due to the relations \(2.1) and \(2.2)
the product $\phi_{\ts kl\ts}(c_k\ts,c_l)\cdot\The$ is divisible on the right
by \(2.5). Since $\Phi=\Up\,\The$ we obtain the first statement of 
Corollary 2.5.
Note that the image of \(2.5) with respect to the antiautomorphism $\io$
is \(2.6). Since $\Phi$ is invariant with respect to $\io$ we get the second 
statement
of Corollary 2.5
\enddemos

\nt
Consider again the rational function $\phik$ appearing at either side of
\(2.1).
The value at $u=w$ of its restriction to $u=v-1$ is not divisible on the 
right by $\phi_{jk}(v,v-1)=1-(j\ts k)$; see proof of Lemma 2.1.
In the proof of Proposition~2.7 we will use the following instant 
observation.

\proclaim{Lemma 2.6}
The value of the function $\phi_{ijk}(v-1,v,w)\,\phi_{kj}(w,v)$ at $w=v-1$
coincides with that of the function $-2\,(ik)\cdot\phi_{kj}(w,v)$. 
\endproclaim

\nt
The next proposition is the central part of the proof of Theorem 2.2\,;
cf. [JKMO].

\proclaim{Proposition 2.7}
Let the numbers $k$ and $k+1$ stand in the same column of the tableau $\Lac$.
Then the element $\Up$ in $\CSn$ is divisible on the left by $1-(k\ts, k+1)$.
\endproclaim

\demo{Proof}
Observe first that Proposition 2.7 follows from its particular case $k+1=n$.
Indeed, let $\nu$ be the shape of the tableau obtained from $\Lac$ by 
removing each of the numbers $k+2\lc n$. Then 
$$ 
\Upz=\Up_\nu(z_1\lc z_{k+1})
\cdot
\prod_{(i,j)}^\rightarrow\ts\ts
\phi_{ij\ts}(c_i+z_i\ts,c_j+z_j)
$$
where $j=k+2\lc n$ and $i$ runs through the sequence $\A_j$. 
Denote by $\Up^{\ts\prime}$ the value at $z_1=\ldots=z_{k+1}$ of the 
restriction of $\Up_\nu(z_1\lc z_{k+1})$ to $\Z$. According to our proof of 
Proposition 2.3 then $\Up=\Up^{\ts\prime}\ts\Up^{\ts\prime\prime}$ for 
certain element $\Up^{\ts\prime\prime}$ in $\CSn$.

{}From now on we will assume that $k+1=n$. Since the element $\The$ is 
invertible, it suffices to prove that $\Up\ts\The=\Phi$ is divisible by 
$1-(n-1\ts, n)$ on the left. But the element $\Phi$ is invariant with respect 
to the antiautomorphism $\io$. We will prove that $\Phi$ is divisible by 
$1-(n-1\ts, n)$ on the right. That is, 
$$
\Phi\cdot\phi_{n\ts,n-1\ts}(c_n\ts,c_{n-1})= 
\Phi\cdot\bigl(1+(n-1\ts,n)\bigr)=0.
\Tag{2.7}
$$

Suppose that the number $n$ appears in the $s$-th column and the $t$-th row
of the tableau $\Lac$. Let $i_1\lc i_s$ be all the numbers in the $t$-th row.
So we have $i_s=n$.
Then due to the relations \(2.1) and \(2.2) for certain element 
$\The^{\ts\prime}$ in $\CSn$ 
$$
\The\cdot\phi_{n\ts,n-1\ts}(c_n\ts,c_{n-1})= \prod_{p<s}^\rightarrow\ts
\phi_{\ts i_p,\ts n-1\ts}(c_{i_p}\ts,c_{n-1}) \cdot
\phi_{n\ts,n-1\ts}(c_n\ts,c_{n-1})
\cdot
\The^{\ts\prime}.
$$
Therefore to get \(2.7)
we have to prove that
$$
\Up\ts\ts\cdot\ts
\prod_{p<s}^\rightarrow\ts
\phi_{\ts i_p,\ts n-1\ts}(c_{i_p}\ts,c_{n-1}) \cdot
\phi_{n\ts,n-1\ts}(c_n\ts,c_{n-1})=0.
\Tag{2.8}
$$
We will prove it by induction on $s$.
If $s=1$ then $\Upz=\Phiz$ so that $\Up=\Phi$. Moreover, 
then none of the pairs $(i,j)$ in \(2.3) is singular. 
Then $\Phi$ has the form
$$
\Phi^{\ts\prime}\cdot\phi_{n-1\ts,n\ts}(c_{n-1}\ts,c_n)= 
\Phi^{\ts\prime}\cdot\bigl(1-(n-1\ts,n)\bigr) 
$$
for certain element $\Phi^{\ts\prime}$ in $\CSn$. 
So we get the equality
$\Up\cdot\phi_{n\ts,n-1\ts}(c_n\ts,c_{n-1})=0$. 

Now suppose that $s>1$. We have to prove that the restriction to $\Z$ 
of the product
$$
\Upz\ts\ts\cdot\ts
\prod_{p\le s}^\rightarrow\ts
\phi_{\ts i_p,\ts n-1\ts}(c_{i_p}+z_{i_p}\ts,c_{n-1}+z_{n-1}) 
%\cdot
%\phi_{n\ts,n-1\ts}(c_n+z_n\ts,c_{n-1}+z_{n-1}) 
\Tag{2.10}
$$
vanishes at $z_1=\ldots=z_n$.
Denote $i_{s-1}=m$. The number $m-1$
appears in the $(s-1)$-\ts th column and the $(t-1)$-\ts th row of $\Lac$.
So we have $c_{m-1}=c_n$.
Let $\mu$ be the shape of the tableau obtained from $\Lac$ by removing each 
of the numbers $m+1\lc n$. Then the function $\Upz$ has the form
$$
\gather
\Up_\mu(z_1\lc z_m)\cdot\Psi(z_1\lc z_{n-1}) \cdot
\phi_{\ts m-1,\ts n-1\ts}(c_{m-1}+z_{m-1}\ts,c_{n-1}+z_{n-1})\times \\
\Chi(z_1\lc z_{n})
\cdot
\phi_{\ts m-1,\ts n\ts}(c_{m-1}+z_{m-1}\ts,c_{n}+z_{n}) \,
\phi_{\ts n-1,\ts n\ts}(c_{n-1}+z_{n-1}\ts,c_{n}+z_{n})\times \\
\prod_{p<s}^\rightarrow\ts
\phi_{\ts i_p n\ts}(c_{i_p}+z_{i_p}\ts,c_{n}+z_{n}). \endgather
$$
Here we have denoted by $\Psi(z_1\lc z_{n-1})$ the product 
$$
\prod_{(i,j)}^\rightarrow\ts\ts
\phi_{ij\ts}(c_i+z_i\ts,c_j+z_j)\ts;
\qquad
j=m+1\lc n-1
\Tag{2.85}
$$
where $i$ runs through $\A_j$ but
$(i,j)\neq(m-1\ts,n-1)$. Further, we have denoted 
$$
\Chi(z_1\lc z_{n})=
\prod_{(i,n)}^\rightarrow\ts\ts
\phi_{in\ts}(c_i+z_i\ts,c_n+z_n)
\Tag{2.9}
$$
where $i$ runs through the sequence $\A_n$ but $i\neq m-1\ts,n-1\ts\lc m$.
In particular, any factor in the product \(2.9) commutes with 
$$
\phi_{\ts m-1,\ts n-1\ts}(c_{m-1}+z_{m-1}\ts,c_{n-1}+z_{n-1}) 
$$
due to \(2.2). Therefore the product \(2.10) takes the form
$$
\gather
\Up_\mu(z_1\lc z_m)
\cdot
\Psi(z_1\lc z_{n-1})
\cdot
\Chi(z_1\lc z_{n})
\,\times
%\Tag{2.11}
\\
\phi_{\ts m-1,\ts n-1\ts,n\ts}(c_{m-1}+z_{m-1}\ts,c_{n-1}+z_{n-1},c_{n}+z_{n})
\,\times
\\
\prod_{p<s}^\rightarrow\ts
\phi_{\ts i_p n\ts}(c_{i_p}+z_{i_p}\ts,c_{n}+z_{n}) \cdot
\prod_{p<s}^\rightarrow\ts
\phi_{\ts i_p,\ts n-1\ts}(c_{i_p}+z_{i_p}\ts,c_{n-1}+z_{n-1}) \,\times
\\
\phi_{n\ts,n-1\ts}(c_n+z_n\ts,c_{n-1}+z_{n-1})=
\\
\Up_\mu(z_1\lc z_m)
\cdot
\Psi(z_1\lc z_{n-1})
\cdot
\Chi(z_1\lc z_{n})
\,\times
\Tag{2.11}
\\
\phi_{\ts m-1,\ts n-1\ts,n\ts}(c_{m-1}+z_{m-1}\ts,c_{n-1}+z_{n-1},c_{n}+z_{n}) 
\,\,
\phi_{n\ts,n-1\ts}(c_n+z_n\ts,c_{n-1}+z_{n-1}) \\
\prod_{p<s}^\rightarrow\ts
\phi_{\ts i_p,\ts n-1\ts}(c_{i_p}+z_{i_p}\ts,c_{n-1}+z_{n-1})
 \cdot\prod_{p<s}^\rightarrow\ts
\phi_{\ts i_p n\ts}(c_{i_p}+z_{i_p}\ts,c_{n}+z_{n})\,. 
\endgather
$$
To get the latter equality we used the relations \(2.1) and \(2.2). 
Restriction to $\Z$
of the product of factors in the first line of \(2.11) is regular at 
$z_1=\ldots=z_n$ according to our proof of Proposition 2.3. Each of the
factors in the third line is also regular at $z_1=\ldots=z_n$. 
Therefore due to Lemma 2.6 the restriction of \(2.11) to $\Z$ has the same 
value at $z_1=\ldots=z_n$
as the restriction to $\Z$ of
$$
\gather
-\,2\cdot
\Up_\mu(z_1\lc z_m)
\cdot
\Psi(z_1\lc z_{n-1})
\cdot
\Chi(z_1\lc z_{n})
\,\times
\Tag{2.12}
\\
(m-1,n)\cdot
\phi_{n\ts,n-1\ts}(c_n+z_n\ts,c_{n-1}+z_{n-1}) \,\times
\\
\prod_{p<s}^\rightarrow\ts
\phi_{\ts i_p,\ts n-1\ts}(c_{i_p}+z_{i_p}\ts,c_{n-1}+z_{n-1}) 
\cdot\prod_{p<s}^\rightarrow\ts
\phi_{\ts i_p n\ts}(c_{i_p}+z_{i_p}\ts,c_{n}+z_{n})= \\
-\,2\cdot
\Up_\mu(z_1\lc z_m)
\cdot
\Psi(z_1\lc z_{n-1})
\cdot
\Chi(z_1\lc z_{n})
\,\times
\\
\prod_{p<s}^\rightarrow\ts
\phi_{\ts i_p,\ts m-1\ts}(c_{i_p}+z_{i_p}\ts,c_{n}+z_{n}) \cdot
\prod_{p<s}^\rightarrow\ts
\phi_{\ts i_p,\ts n-1\ts}(c_{i_p}+z_{i_p}\ts,c_{n-1}+z_{n-1}) \,\times
\\
(m-1,n)\cdot \phi_{n,n-1\ts}(c_n+z_n\ts,c_{n-1}+z_{n-1})\,. \endgather
$$
Here each of the factors
$\phi_{\ts i_p,\ts m-1\ts}(c_{i_p}+z_{i_p}\ts,c_{n}+z_{n})$ commutes with
$\Chi(z_1\lc z_{n})$ by the relations \(2.2). 
In each of these factors we can replace $c_{n}+z_{n}$ 
by $c_{m-1}+z_{m-1}$ without affecting the value at $z_1=\ldots=z_n$ 
of the restriction to $\Z$ of \(2.12).
Denote
$$
\Gaz\ts=\ts
\prod_{p<s}^\rightarrow\ts
\phi_{\ts i_p,\ts m-1\ts}(c_{i_p}+z_{i_p}\ts,c_{m-1}+z_{m-1}). $$
According to our proof of Proposition 2.3 it now suffices to 
demonstrate vanishing at $z_1=\ldots=z_{n-1}$ of the 
restriction to $\Z$ of the product 
$$
\Up_\mu(z_1\lc z_m)
\cdot
\Psi(z_1\lc z_{n-1})
\cdot
\Gaz
\ts.
\Tag{2.13}
$$

Consider the product \(2.85). Here the factors corresponding to the pairs
$(i,j)$ are arranged with respect to ordering chosen in the proof of 
Proposition 2.3. Let us now reorder the pairs $(i,j)$ in \(2.85) as follows.
For every $j=m+1\lc m+\lap_{s-1}-t$ change the sequence 
$$
(m-1,j)\ts,(i_1,j)\lc(i_{\ts s-1},j)
\nopagebreak
$$ 
to
$$
(i_1,j)\lc(i_{\ts s-1},j)\ts,(m-1,j)\ts. 
$$ 
Denote by $\Psi^{\ts\prime}(z_1\lc z_{n-1})$ the resulting ordered product.
Then by \(2.1) and \(2.2)
$$
\Psi(z_1\lc z_{n-1})\cdot\Gaz=
\Gaz\cdot\Psi^{\ts\prime}(z_1\lc z_{n-1}). \Tag{2.135}
$$

Now let $(i,j)$ be any singular pair in \(2.85). Let $(h,j)$ be the pair 
following $(i,j)$ in the ordering from the proof of Proposition 2.3.
Then $(h,j)$ follows $(i,j)$ in our new ordering as well. 
Furthermore, due to the relations \(2.1) and \(2.2)
the product
$$
\Up_\mu(z_1\lc z_m)
\cdot
\Gaz
\Tag{2.14}
$$
is divisible on the right by $\phi_{ih}(c_i+z_i,\ts c_h+z_h)$. Therefore
$$
\align
&
\Up_\mu(z_1\lc z_m)
\cdot
\Gaz
\cdot
\Psi^{\ts\prime}(z_1\lc z_{n-1})=
\Tag{2.15}
\\
&
\Up_\mu(z_1\lc z_m)
\cdot
\Gaz
\cdot
\Psi^{\ts\prime\prime}(z_1\lc z_{n-1})
\endalign
$$
where $\Psi^{\ts\prime\prime}(z_1\lc z_{n-1})$ is a rational function which
restriction to is $\Z$ regular at $z_1=\ldots=z_{n-1}$. But by the inductive
assumption the restriction of \(2.14) to $\Z$ vanishes at $z_1=\ldots=z_{m}$.
Thus due to \(2.135) and \(2.15) the restriction of \(2.13) to 
$\Z$ vanishes at $z_1=\ldots=z_{n-1}$ 
\enddemos

\nt
We have shown that the element $\Phi=\Up\,\The$ of $\CSn$ is 
invariant with respect to~the antiautomorphism $\io$. 

\proclaim{Corollary 2.8}
Let the numbers $k$ and $k+1$ stand in the same column of the tableau $\Lac$.
Then the element $\Phi$ in $\CSn$ is divisible on the 
left and on the~right~by
$\phi_{k,\ts k+1}(c_k,c_{k+1})$.
\endproclaim

\nt
The next proposition completes the proof of Theorem 2.2. 

\proclaim{Proposition 2.9}
We have the equality $\Phi=\Phi_\la$.
\endproclaim

\demo{Proof}
Due to Propositions 2.4 and 2.7 we have the equality 
$\Up=Q_\la\ts X_\la\ts P_\la$ for some element $X_\la\in\CSn$. 
It is a classical fact that one can assume here $X_\la\in\CC$; 
see for instance [GM]. Due to Corollary 2.8 the decomposition 
$\Phi=\Up\,\The$ now implies that $\Phi$ equals $Q_\la\ts P_\la\ts Q_\la$
up to a multiple from $\CC$. We will show that in the expansion of 
$\Phi\in\CSn$ with respect to the basis of $g\in S_n$ 
the coefficient at identity is $1$.

Let $g_0$ be the element of the maximal length in $S_n$. 
Consider the product \(2.3) with the initial lexicographical ordering 
of the pairs $(i,j)$. For each $k=1\lc n-1$ denote 
$$
\phi_{k}(u,v)=\phi_{k,k+1}(u,v)\cdot(k,\ts k+1)=(k,\ts k+1)-\frac1{u-v}. 
$$
Then $$
\Phiz\,g_0=
\prod_{(i,j)}^\rightarrow\ts\ts
\phi_{j-i\ts}(c_i+z_i\ts,c_j+z_j).
$$
But $$
g_0=\prod_{(i,j)}^\rightarrow\ts\ts (j-i,j-i+1) $$
is a reduced decomposition in $S_n$.
Therefore in the expansion of the element $$
\Phiz\,g_0\in\ts\CC(z_1\lc z_n)\!\cdot\! S_n $$
with respect to the basis of $g\in S_n$ the coefficient at $g_0$ is $1$.
So is the coefficient at $g_0$ in the expansion of the element 
$\Phi\,g_0\in\CSn$
\enddemos

\nt Continuation of \(2.3) to $z_1=\ldots=z_n$ along the set $\Z$ 
is called {\it fusion procedure}. In the course of the proof of 
Proposition 2.9 we have established
the following fact.

\proclaim{Corollary 2.10}
We have the equality $\Up=Q_\la\ts P_\la$. \endproclaim

\nt
We will now make a concluding remark about the element $\Up\in\CSn$. 
Consider the {\it row tableau} of the shape $\la$ obtained by filling 
the boxes of $\la$ with $1\lc n$ by rows downwards, 
from the left to the right in every row. Denote by $\Lar$ this tableau. 
It is obtained from the tableau $\Lac$ by a certain permutation of the
numbers $1\lc n$. Denote by $\si_\la$ this permutation. 
For the vectors $u_\Lac,u_\Lar$ of the Young basis in $U_\la$
we have
$\langle\ts\si_\la\ts u_\Lac\ts,u_\Lar\ts\rangle\neq0$ due to [Y2]. 
The element
$$
\Up_\la=\sum_{\si\in S_n}\ts
\frac
{\langle\ts\si\ts u_\Lac\ts,u_\Lar\ts\rangle} 
{\langle\ts\si_\la\ts u_\Lac\ts,u_\Lar\ts\rangle} \,\si\1\si_\la
$$
of the group ring $\CSn$ does not depend on the choice of the vectors 
$u_\Lac,u_\Lar$. In fact one has $\Up_\la=Q_\la\ts P_\la$; 
see for instance [GM]. Thus $\Up=\Up_\la$ by Corollary 2.10. 

\proclaim{Lemma 2.11}
\hskip-2pt
Restriction of $\phik\hskip-1pt$ to the set of all $(u,v,w)\hskip-2pt$
with \hbox{$w=v\pm1\ts,\hskip-1pt$} is regular at $u=w$. 
\endproclaim

\demo{Proof}
As a direct calculation shows, the restriction 
$$
\phi_{ijk}(u,v,v\pm1)=
\left(1-\frac{(i\ts j)+(i\ts k)}{u-v}\right) \!\cdot\!
\Bigl(1\pm(j\ts k)\Bigr).
\nopagebreak
$$
The latter rational function of $u,v$ is manifestly regular at $u=v\pm1$
\enddemos

\nt
We will end up this section with one more proposition. 
It will be used in the next section. 
Let the superscript ${}^\vee$ denote the embedding of the group 
$S_n$ to $S_{n+1}$ determined by assigning to the transposition of 
$i$ and $j$ that of $i+1$ and $j+1$. 

\proclaim{Proposition 2.12}
We have equality of rational functions in $u$ valued in 
$\CC\!\cdot\!S_{n+1}$
$$
\prod^\rightarrow_{1\le i\le n}\ \phi_{1,i+1}(u\ts,c_i) \cdot
\Phi_\la^\vee=
\biggl(1-\sum_{1\le i\le n}\frac{(1,i+1)}u\,\biggr) \cdot
\Phi_\la^\vee
\ts.
\Tag{2.99}
$$
\endproclaim

\demo{Proof}
Denote by $\Xiu$ the rational function at the left hand side of \(2.99).
The value of this function at $u=\infty$ is $\Phi_\la^\vee$. Moreover, 
the residue of $\Xiu$ at $u=0$ is $$
-\sum_{1\le i\le n}{(1,i+1)}\cdot\Phi_\la^\vee. 
$$
It remains to prove that $\Xiu$
has a pole only at $u=0$ and this pole is simple. Let an index 
$i\in\{2\lc n\}$ be fixed. The factor $\phi_{1,i+1}(u\ts,c_i)$
in \(2.99) has a pole at $u=c_i$. We shall prove that when estimating from
above the order of the pole at $u=c_i$ of $\Xiu,$ that factor does not count.

Suppose that the number $i$ does not appear in the first row of the tableau 
$\Lac$. Then the number $i-1$ appears in $\Lac$ straight above $i$. 
In particular, $c_{i-1}=c_i+1$. By Corollary 2.8 the element $\Phi_\la^\vee$ 
is divisible on the left by $\phi_{i,i+1}(c_{i-1}\ts,c_i)$.
But the product
$$
\phi_{\ts 1i}(u\ts,c_{i-1})\ts\ts
\phi_{\ts 1,i+1}(u\ts,c_{i})\ts\ts
\phi_{i,i+1}(c_{i-1}\ts,c_i)
$$
is regular at $u=c_i$ due to Lemma 2.11. 

Now suppose that the number $i$ appears in the first row of $\Lac$. 
Let the number $k$ be next to the left of $i$ in the first row of $\Lac$.
Then $c_k=c_i-1$.
By Corollary~2.5 the element $\Phi_\la^\vee$ is divisible on the left by
$$
\prod_{k<j<i}^\leftarrow\ts
\phi_{\ts j+1,i+1\ts}(c_j\ts,c_i)
\cdot
\phi_{\ts k+1,i+1\ts}(c_k\ts,c_i).
\Tag{2.98}
$$
Consider the product of factors in \(2.99) $$
\phi_{\ts 1,k+1\ts}(u\ts,c_k)
\cdot\!\!
\prod_{k<j<i}^\rightarrow\!
\phi_{\ts 1,j+1\ts}(u\ts,c_j)
\cdot
\phi_{\ts 1,i+1\ts}(u\ts,c_i).
\Tag{2.97} $$
Multiplying the product \(2.97) on the right by \(2.98) and using 
\(2.1),\(2.2) we get
$$
\gather
\prod_{k<j<i}^\leftarrow\ts
\phi_{\ts j+1,i+1\ts}(c_j\ts,c_i)
\ \times
\\
\phi_{\ts 1,k+1}(u,c_{k}u)\ts\ts
\phi_{\ts 1,i+1}(u,c_{i})\ts\ts
\phi_{k+1,i+1}(c_{k}\ts,c_i)
\\
\times\ \prod_{k<j<i}^\rightarrow\!
\phi_{\ts 1,j+1\ts}(u\ts,c_j).
\endgather
$$
The latter product is regular at $u=c_i$ by Lemma 2.11. 
The proof is complete
\enddemos

\newpage

\centerline{\S\ts\bf3.\ Yangian of the general linear Lie algebra}\section{\,}
\kern-20pt

\nt
In this section we will collect several known facts from [C2\ts,D1] about the
{\it Yangian}
of the Lie algebra $\glN$. 
This is a complex associative unital algebra $\YN$ with the
countable
set of generators $T^{(s)}_{ij}$ where $s=1,2,\ts\ldots$ and 
$i,j=1,\ts\dots\ts,N$. The de\-fin\-ing relations in the algebra $\YN$ are
$$
[\ts T_{ij}^{(r+1)},T_{kl}^{(s)}\ts]-
[\ts T_{ij}^{(r)},T_{kl}^{(s+1)}\ts]=
T_{kj}^{(r)}\ts T_{il}^{(s)}-T_{kj}^{(s)}\ts T_{il}^{(r)};
\quad r,s=0,1,2,\ldots\ts
\Tag{3.1}
$$
where $T_{ij}^{(0)}=\de_{ij}\cdot1$.
We will also
use the following matrix form of these relations. 
\Par
Let $E_{ij}\in\EndCN$ be the standard matrix units. Introduce the element
$$
P=\sum_{i,j}\ts E_{ij}\ot E_{ji}\in\EndCN\ot\EndCN\ts $$
where the indices $i,j$ run through $1\lc N$. Consider the 
{\it Yang $R$-matrix\,},
it is the rational function of two complex variables $u\,,v$ valued in 
$\EndCN\ot\EndCN$
$$
R\ts(u\,,v)=\id-\frac P{u-v}\ts.
$$
Introduce the formal power series in~$u^{-1}$ $$
T_{ij}(u)=
T_{ij}^{(0)}+T_{ij}^{(1)}\ts u^{-1}+T_{ij}^{(2)}\ts u^{-2}+\ldots $$
and combine all these series into the single element of 
$\EndCN\ot\YN\ts[[u^{-1}]]$
$$
T(u)=\sum_{i,j}\ts E_{ij}\ot T_{ij}(u).
$$

For any associative unital algebra $\operatorname{A}$ 
denote by $\iota_s$ its embedding into the tensor product 
$\operatorname{A}^{\!\ot n}$
as the $s$-\ts th tensor factor:
$$
\iota_s(X)=1^{\ot\ts (s-1)}\ot X\ot1^{\ot\ts(n-s)}\,, \quad
X\in\operatorname{A}\,;
\qquad
s=1,\dots,n.
$$
We will also use various embeddings of the algebra 
$\operatorname{A}^{\!\ot\ts m}$ into $\operatorname{A}^{\!\ot\ts n}$ 
for any $m\leqslant n$.
For $s_1\lc s_m\in\{1\lc n\}$ pairwise distinct and 
$X\in\operatorname{A}^{\!\ot m}$
we put
$$
X_{s_1\ldots s_m}=
\iota_{s_1}\!\ot\ldots\ot\iota_{s_m}(X)\in\operatorname{A}^{\!\ot n}. 
$$
Denote
$$
T_s(u)=\iota_s\ot\id\bigl(T(u)\bigr)\in\EndCN^{\ot n}\ot\YN\,[[u\1]]. 
$$
In this notation the defining relations \(3.1) can be rewritten as the
single equation
$$
R(u,v)\ot1\cdot\ts T_1(u)\ts T_2(v)=
T_2(v)\ts T_1(u)\cdot R(u,v)\ot1.
\nopagebreak
\Tag{3.2}
$$
After multiplying each side of \(3.2) by $u-v$ it becomes a relation in
the algebra
$$
\EndCN\ot\EndCN\ot\YN\,((u\1,v\1)).
\nopagebreak
$$
For further comments on the definition of the algebra $\YN$ see [MNO].

The relation \(3.2) implies that for any $z\in\CC$
the assignment $T_{ij}(u)\mapsto T_{ij}(u-z)$ determines an automorphism of 
the algebra $\YN$. Here the formal series in $(u-z)\1$ should be re-expanded 
in $u\1$. We will denote by $\tau_z$ this automorphism. 

We will also regard $E_{ij}$ as generators of the universal enveloping 
algebra $\UN$.
The algebra $\YN$ contains $\UN$ as
a subalgebra: due to \(3.1) the assignment $E_{ij}\mapsto -\,T_{ji}^{(1)}$ 
defines the embedding. Moreover, there is a homomorphism
$$
\pi:
\ts\YN\to\UN:\ts T_{ij}(u)\mapsto \de_{ij}-E_{ji}\ts u^{-1}. 
\Tag{3.5}
$$
This homomorphism is by definition
identical on the subalgebra $\UN$. It is called the 
{\it evaluation homomorphism} for the algebra $\YN$. 
We will regard $\CC^N$ as $\YN$-module by virtue of this homomorphism. 
Denote by $V(z)$ the $\YN$-module obtained from $\CC^N$ by applying the 
automorphism $\tau_z$.
The action of the generators $T_{ij}^{(s)}$ in $V(z)$ can be determined by 
the assignment
$$
\EndCN\ot\YN\ts[[u^{-1}]]
\rightarrow
\EndCN\ot\EndCN\ts[[u^{-1}]]
\,:\
T(u)\mapsto R(u,z).
$$

Furthermore,
there is a natural Hopf algebra structure on $\YN$. Again due to \(3.2) the 
comultiplication
$\Delta:\YN\to\YN\ot\YN$ can be defined by \smallskip
$$
T_{ij}(u)\mapsto\sum_{h}\ts
T_{ih}(u)\ot T_{hj}(u).
\Tag{3.3}
$$
Here the tensor product is taken over the subalgebra $\CC[[u^{-1}]]$ in 
$\YN\ts[[u\1]]$ and $h$ runs through $1\lc N$. 
For any $z_1\lc z_n\in\CC$ consider
the $\YN$-module
$V(z_1)\ot\ldots\ot V(z_n)$.
By the definition \(3.3)
the action of the generators $T_{ij}^{(s)}$ in the space $(\CC^N)^{\ot n}$ 
of this module
can be determined by the assignment
$$
\gather
\EndCN\ot\YN\ts[[u^{-1}]]
\ \rightarrow\
\EndCN\ot\EndCN^{\ot n}\ts[[u^{-1}]]\ts:\ \ \Tag{3.6}
\\
T(u)
\,\mapsto\,
R_{12}(u,z_1)\ts\ldots\ts R_{1,n+1}(u,z_n). \endgather
$$

\Par
The symmetric group $S_n$ acts in the space $(\CC^N)^{\ot n}$ 
by permutations of the tensor factors:
$$
(k,l)\ \mapsto\ P_{\ts kl}\ts\in\,\EndCN^{\ot n}\ ; \qquad
k\neq l\,.
$$
The function $R_{kl}(u,v)$ valued in $\EndCN^{\ot n}$ 
corresponds to the function $\phi_{kl}(u,v)$
%from the previous section
valued in the group ring $\CSn$\ts.
Let $\la$ be a Young diagram with $n$ boxes and not more than $N$ rows.
Denote by $F_\la$ the element of $\EndCN^{\ot n}$ corresponding to
$$
\frac{\dim U_\la}{n!}\,\ts\Phi_\la\in\CSn\ts. 
$$
So $F_\la^{\ts2}=F_\la$
and the image of $F_\la$ in $(\CC^N)^{\ot n}$ is an irreducible 
$\glN$-submodule~[Y1]. Denote this image by $V_\la$. We will identify 
the algebra $\End(V_\la)$ with the subalgebra in $\EndCN^{\ot n}$ which 
consists of all the elements of the form $F_\la\ts X\ts F_\la$.

By virtue of the homomorphism \(3.5) we can regard $V_\la$ as an 
$\YN$-module. Denote the action of $\YN$ in the module $V_\la$ by 
$\pi_\la$. 
As well as in the previous section, 
let $c_k$ be the content of the box with the number $k$ in the column 
tableaux $\Lac$. The next proposition is contained in [C2]. 

\proclaim{Proposition 3.1}
The action $\pi_\la$ of $\YN$
can be described as the restriction
to the vector subspace $V_\la$ in $(\CC^N)^{\ot n}$ of the action \(3.6) 
with
$z_1=c_1\ts\lc z_n=c_n$.
\endproclaim

\demo{Proof}
Consider the vector space $(\CC^N)^{\ot n}$ as an $\YN$-module by
virtue of 
\(3.5).
The action of $\YN$ in the latter module can be determined by the 
assignment
$$
\gather
\EndCN\ot\YN\ts[[u^{-1}]]
\ \rightarrow\
\EndCN\ot\EndCN^{\ot n}\,[[u^{-1}]]\ts:\ \ \Tag{3.7}
\\
T(u)
\,\mapsto\,
1-\sum_{1\le k\le n}\frac{P_{1,k+1}}u
\endgather
$$
The action of $\YN$ in $V_\la$ is by definition the restriction of \(3.7) 
to
the subspace $V_\la$ in $(\CC^N)^{\ot n}$. So it suffices to prove the 
equality of rational functions in $u$ $$
\Bigl(1-\sum_{1\le k\le n}\frac{P_{1,k+1}}u\,\Bigr) \cdot
\bigl(\id\ot F_\la\bigr)
\ =
\prod_{1\le k\le n}^\rightarrow
R_{1,k+1}(u,c_k)
\cdot
\bigl(\id\ot F_\la\bigr).
$$
valued in $\EndCN^{\ot(n+1)}$. But this equality is afforded by 
Proposition 
2.12
\enddemos

\nt
Denote by $F_N$ the antisymmetrization map in $\EndCN^{\ot N}$. 
Thus $F_N=F_\la$ for $n=N$ and the diagram $\la$ consisting of one column 
only. Consider the element
$$
F_N\ot1\cdot T_1(u)\ts\ldots\ts T_N(u-N+1) \Tag{3.9}
$$
of the algebra $\EndCN^{\ot N}\ot\YN[[u\1]]$. 
Due to Theorem 2.2 and the relations \(3.2) this element is divisible by 
$F_N\ot1$ also on the right. Since the map $F_N$ is one-dimensional, the 
element \(3.9) has the form $F_N\ot D(u)$ for a certain series 
$D(u)\in\YN[[u\1]]$. This series is called the {\it quantum determinant} 
for $\YN$. The coeficients $D_1,D_2,\ldots$
of this series at $u\1, u^{-2}\ldots$
are free generators of the centre of the algebra $\YN;$ see [MNO] for the 
proof. The definitions \(3.3) and \(3.9) imply
$$
\Delta\bigl(D(u)\bigr)=D(u)\ot D(u).
$$
Hence the centre of the algebra $\YN$ is a Hopf subalgebra. 

The relations \(3.2) imply that
for any formal series $f(u)\in 1+u^{-1}\ts\CC[[u^{-1}]]$ the assignment
$T_{ij}(u)\mapsto f(u)\cdot T_{ij}(u)$
determines an automorphism of the algebra $\YN$. We will denote by $\om_f$ 
this automorphism.
Consider the fixed point subalgebra in $\YN$ with respect to all 
the automorphisms $\om_f$. This subalgebra is called the {\it Yangian} 
of the simple Lie algebra $\slN$ and denoted by $\YslN\ts;$ cf.~[D1]. 

\proclaim{Proposition 3.2}
{\hskip-1pt}The algebra $\YslN\!$ is a Hopf subalgebra in $\hskip-1pt\YN.$ 
\hskip-1pt\text{The algebra} $\YN$ is isomorphic to the tensor product of 
its centre and its subalgebra $\YslN.$ \endproclaim

\demo{Proof}
We will follow the argument from [MNO]. Observe that by our definition 
$$
\om_f\bigl(D(u)\bigr)=D(u)\cdot f(u)\ts\ldots\ts f(u-N+1). 
$$
Determine the series
$A(u)\in\YN[[u\1]]$ by the equation
$$
D(u)=A(u)\ts\ldots\ts A(u-N+1).
$$
Then we have $\om_f\bigl(A(u)\bigr)=f(u)\ts A(u).$ 
The coefficients $A_1, A_2\ldots$
of the series $A(u)$ at $u\1, u^{-2}\ldots$ are free generators of the 
centre of $\YN.$ Moreover, we have
$$
\Delta\bigl(A(u)\bigr)=A(u)\ot A(u).
\Tag{3.222}
$$

\newpage

Every coefficient of the series $T_{ij}(u)\ts A(u)\1$ belongs to the 
subalgebra $\YslN.$ Therefore the algebra $\YN$ 
is generated by its centre and its subalgebra $\YslN.$
Now suppose that for some
positive integer $s$
there exists a non-zero polynomial $\Cal Q$ in $s$ variables with 
the coefficients from $\YslN$ such that $\Cal Q(A_1\lc A_s)=0$. 
Choose the number $s$ to be minimal.
But for $f(u)=1+a\ts u^{-s}$
with any $a\in\CC$ we have
$$
\om_f:\,\Cal Q(A_1\lc A_s)\,\mapsto\,\Cal Q(A_1\lc A_s+a). $$
Thus $\Cal Q(A_1\lc A_s+a)=0$ for all $a\in\CC$. So we may assume that 
the polynomial $\Cal Q$ does not depend on the last variable, 
which contradicts to the choice of $s$. 
This contradiction completes the proof of the second statement of 
Proposition 3.2. 

In particular, the coefficients of all the series $T_{ij}(u)\ts A(u)\1$ 
generate the algebra $\YslN$. Now the first statement of 
Proposition 3.2 follows from \(3.3) and \(3.222)
\enddemos

\proclaim{Corollary 3.3}
The centre of the algebra $\YslN$ is trivial. \endproclaim

\nt
Note that the subalgebra $\YslN$ in $\YN$ is preserved by the 
automorphism $\tau_z\ts$ for any $z\in\CC$. Denote by $\De^\prime$ 
the composition of $\De$ with the permutation of tensor factors in 
$\YN\ot\YN$.
According to [D1]
there is a unique formal series
$$
\SSS(z)\ts\in\,1+\YslN\ot\YslN\,[[z\1]]\cdot z\1 $$
such that in the algebra $\YslN^{\ot3}\,[[z\1]]$ we have $$
\align
\id\ot\De\ts\bigr(\SSS(z)\bigl)&=\SSS_{\ts12}(z)\,\SSS_{\ts13}(z), 
\\
\De\ot\id\ts\bigr(\SSS(z)\bigl)&=\SSS_{\ts13}(z)\,\SSS_{\ts23}(z) 
\Tag{3.333}
\endalign
$$
and for any element $Y\in\YslN$ we also have $$
\SSS(z)\cdot\id\ot\tau_z\bigl(\De^\prime(Y)\bigr) 
=\id\ot\tau_z\bigl(\De(Y)\bigr)\cdot\SSS(z)\ts. 
\Tag{3.4}
$$
The series $\SSS(z)$ is called the {\it universal $R$-matrix\,} 
for the Hopf algebra $\YslN$.
For any $v\in\CC$ denote by $\pi_v$ the action of the algebra $\YN$ 
in the module $V(v)$.

\proclaim{Lemma 3.4}
We have the equality
$\pi_v\ot\id\ts\bigl(\SSS(z)\bigr)=T(v-z)\cdot\ts\id\ot A(v-z)\1\ts.$ 
\endproclaim

\demo{Proof}
Let us keep the parameter $v\in\CC$ fixed. Denote respectively by 
$T^{\ts\prime}(z)$ and $T^{\ts\prime\prime}(z)$ the formal series in 
$z\1$ at the left and the right hand side of the equality to prove.
Due to \(3.4) and Proposition~3.2
the series $T^{\ts\prime}(z)$ satisfies the equation $$
T^{\ts\prime}(z)\cdot\pi_v\ot\tau_z\bigl(\De^\prime(Y)\bigr) =
\pi_v\ot\tau_z\bigl(\De(Y)\bigr)\cdot T^{\ts\prime}(z) \Tag{3.555}
$$
for any element $Y\in\YN$.
We will assume that here $Y$ is one of the generators $T_{ij}^{(s)}$. 
By \(3.3)
we have the equalities in $\EndCN\ot\EndCN\ot\YN[[u\1]]$ $$
\gather
\id\ot(\pi_v\ot\tau_z\ts\crc\ts\De\phantom{^\prime})\bigl(T(u)\bigr)= 
R_{12}(u,v)\ot1
\cdot
T_1(u-z),
\\
\id\ot(\pi_v\ot\tau_z\ts\crc\ts\De^\prime)\bigl(T(u)\bigr)= \,T_1(u-z)
\cdot
R_{12}(u,v)\ot1.
\endgather
\nopagebreak
$$
So the collection of all equations \(3.555) is equivalent to the single 
equation for $T^{\ts\prime}(z)$ 
$$
R_{12}(u+z,v)\ot1
\cdot
T_1(u)
\,
T^{\ts\prime}_2(z)
=
T^{\ts\prime}_2(z)
\ts
\,T_1(u)
\cdot
R_{12}(u+z,v)\ot1\ts,
\Tag{3.5555}
$$
which after being multiplied by $u+z-v$ becomes an equation in the algebra
$$
\EndCN\ot\EndCN\ot\YN\,((u\1,z\1)).
$$
But the element $T^{\ts\prime\prime}(z)$ satisfies the same equation. 
Indeed,
since the coefficients of the series $A(v-z)$ in $z\1$ are central in 
$\YN$ we have due to \(3.2)
$$
R_{12}(u+z,v)\ot1
\cdot
T_1(u)
\,
T^{\ts\prime\prime}_2(z)
=
T^{\ts\prime\prime}_2(z)
\ts
\,T_1(u)
\cdot
R_{12}(u+z,v)\ot1\ts.
\Tag{3.55555}
$$
Consider the series
$$
X(z)=T^{\ts\prime}(z+v)\ts\,T^{\ts\prime\prime}(z+v)\1
$$
in $z\1\!$ with the coefficients in $\EndCN\ot\YslN.$ 
\!By comparing \(3.5555),\(3.55555)~we~get 
$$
\bigl[\ts
R_{12}(u+z+v,v)\ot1
\cdot
T_1(u)
\ts,\ts
X_2(z)
\ts\bigr]=0\ts.
\Tag{3.555555}
$$

We will write
$$
X(z)=\sum_{i,j}\ts E_{ij}\ot X_{ij}(z)
$$
where
$$
X_{ij}(z)=
X_{ij}^{(0)}+X_{ij}^{(1)}\ts z^{-1}+X_{ij}^{(2)}\ts z^{-2}+\ldots 
$$
for some $X_{ij}^{(s)}\in\YslN$.
Let us multiply the equation \(3.555555) by $u+z$. 
Then by considering the coefficient at $u^0\ts z^{-s}$ we get
$$
\bigl[\ts
T_{ij}^{(1)}\ts,\ts X_{kl}^{(s)}
\ts\bigr]=
\de_{kj}\,X_{il}^{(s)}-\de_{il}\,X_{kj}^{(s)}\ts;\quad s=1,2,\ldots\,\ts.
\Tag{3.5555555}
$$
Further, by considering the coefficient at $u^{-1}\ts z^{-s}$ we get 
$$
\bigl[\ts
T_{ij}^{(2)}\ts,\ts X_{kl}^{(s)}
\ts\bigr]+
\bigl[\ts
T_{ij}^{(1)}\ts,\ts X_{kl}^{(s+1)}
\ts\bigr]=\ts
T_{kj}^{(1)}\ts X_{il}^{(s)}-\ts
\de_{il}\cdot\sum_h\ts X_{kh}^{(s)}\,T_{hj}^{(1)} \Tag{3.55555555}
$$
where the index $h$ runs through $1\lc N$. 

Let us now prove by induction on $s=0,1,2,\ldots$ that 
$X_{kl}^{(s)}=\de_{kl}\,X_s$ for some element $X_s\in\YslN$ which 
commutes with every $T_{ij}^{(1)}\ts$. The condition $$
X(z)\in1+\EndCN\ot\YslN\ts[[z\1]]\cdot z\1 $$
provides the base for our induction.
Let us consider the equation \(3.55555555). 
Due to \(3.5555555) and to the inductive assumption it takes the form 
$$
\de_{kl}\,\bigl[\ts
T_{ij}^{(2)}\ts,\ts X_s
\ts\bigr]+
\de_{kj}\,X_{il}^{(s+1)}-\de_{il}\,X_{kj}^{(s+1)} =0\ts.
\nopagebreak
$$
The latter equation shows that for $k\ne l$ we have $
X_{kk}^{(s+1)}\!=X_{ll}^{(s+1)}
$
and $X_{kl}^{(s+1)}\!=0.$
Thus $X_{kl}^{(s+1)}=\de_{kl}\,X_{s+1}$ for some $X_{s+1}\in\YslN$. 
Now by using \(3.5555555) with $s+1$ instead of $s$ we get the 
equality $[\ts T_{ij}^{(1)}\ts,\ts X_{s+1}\ts]=0$. Therefore
$$
X(z)=\id\ot\bigl(\ts
1+X_1\ts z^{-1}+X_2\ts z^{-2}+\ldots
\ts\bigr)
\nopagebreak
$$
for some elements $X_1\ts,X_2\ts,\ldots\in\YslN.$ 
These elements are central due to \(3.555555). 
\line{But the centre of the algebra $\YslN$ is $\CC$ by Corollary 3.3. 
So $\!X_1\ts,\!X_2\ts,\ldots\in\CC.$}

\newpage

Thus the series $T^{\ts\prime}(z)$ and $T^{\ts\prime\prime}(z)$ 
coincide up to a multiple from $\CC[[z\1]]$. 
That multiple must be $1$ since the equation
$$
\pi_v\ot\De\ts\bigr(\SSS(z)\bigl)=
\pi_v\ot\id\ot\id\ts\bigl(\SSS_{\ts12}(z)\ts\SSS_{\ts13}(z)\bigr) 
\nopagebreak
$$
for $\pi_v\ot\id\,\bigr(\SSS(z)\bigl)=T^{\ts\prime}(z)$ is also satisfied 
for $T^{\ts\prime\prime}(z)$.
The latter fact follows from the definition \(3.3) and from \(3.222) 
\enddemos

\nt
Now introduce the formal series in $z\1$ with the coefficients in the 
centre of $\YN$
$$
D_\la(z)=A(c_1-z)\ts\ldots\ts A(c_n-z).
\Tag{3.99}
$$
Introduce also the element of the algebra $\EndCN^{\ot n}\ot\YN[[z\1]]$ 
$$
T_\la(z)=F_\la\ot1\cdot T_1(c_1-z)\ldots T_n(c_n-z). \Tag{3.8}
$$
By Theorem 2.2 and the relations \(3.2)
the latter element is divisible by $F_\la\ot1$ also on the right. 
Thus we can regard $T_\la(z)$ as an element of $\End(V_\la)\ot\YN[[z\1]]$.
By \(3.333) and by Lemma 3.4 we have the following corollary to 
Proposition 3.1.

\proclaim{Corollary 3.5}
We have the equality
$\pi_\la\ot\id\,\bigr(\SSS(z)\bigl)=T_\la(z)\cdot\ts\id\ot D_\la(z)\1$. 
\endproclaim

\nt
Next theorem is the main result of the present section. 

\proclaim{Theorem 3.6}
There is a series $\R(z)\in1+\YN\ot\YN\ts[[z\1]]\cdot z\1$ 
satisfying \(1.0) such that the equality 
$\pi_\la\ot\id\,\bigr(\R(z)\bigl)=T_\la(z)$ holds for any 
Young diagram~$\la$.
\endproclaim

\demo{Proof}
We will prove that there is a formal power series 
$\Cal D(z)$ in $z\1$~with~coefficients
in the tensor square of the centre of $\YN$ and leading 
term $1$ such~that
$\pi_\la\ot\id\,\bigr(\Cal D(z)\bigl)=\id\ot D_\la(z)$ 
for any Young diagram~$\la$. Then due to Proposition~3.2 and 
Corollary 3.5 it will suffice to put $\R(z)=\SSS(z)\,\Cal D(z).$ 

Let each of the indices $d_1,d_2,\ldots$ run through 
$0,1,2,\ldots$ but only finite number of them differ from zero. 
Put $d_1+d_2+\ldots=d$.
Introduce the formal series in $z\1$
$$
f_{d_1,d_2,\ldots}(z)=
\sum_{s_1\lc s_n}\ts
(c_1-z)^{-s_1}\ldots\,(c_n-z)^{-s_n}
\Tag{3.999}
$$
where the sum is taken over all the sequences $s_1\lc s_n$ of the 
numbers $0,1,2,\ldots$ such that the multiplicities of $1,2,\ldots$ are 
$d_1,d_2,\ldots$ respectively.
Then by \(3.99)
$$
D_\la(z)=\sum_{d_1,\ts d_2,\ts\ldots}
\
f_{d_1,d_2,\ldots}(z)
\cdot
A_1^{d_1}A_2^{d_2}\ldots
\Tag{3.9999}
$$
where the sum is taken over all sequences 
$d_1,d_2,\ldots$ such that $d\le n$.

Now let the indices $d_1,d_2,\ldots$ be fixed but the number $n$ and
the diagram $\la$ vary. For each $i=1\lc N$ denote by $\la_i$
the length of $i$-th row of the diagram~$\la$. 
Let $f_{\la}^{(s)}$ be
the coefficient of the series \(3.999) at~$z^{-s}$. 
It depends on
$\la$ as a fixed polynomial of
$$
c_1^{\ts r}+\ldots+c_n^{\ts r}\,\ts;\qquad r=0,1,2,\ldots\,\ts. $$
So $f_{\la}^{(s)}$ can be expressed as a fixed symmetric polynomial in 
$\la_1-1\lc\la_N-N$.
Therefore it is the eigenvalue in the $\glN$-module $V_\la$ of a 
certain element $F_s$ of the~centre $\ZN$ of the algebra $\UN$ 
which does not depend on the diagram~$\la$. 

\newpage

Note that by \(3.999) we get $f_{\la}^{(s)}=0$ if $s<d_1+2\ts d_2+\ldots$ 
whatever~$\la$~is. Hence for any $s<d_1+2\ts d_2+\ldots$ we have $F_s=0$. 
For $d_1=d_2=\ldots=0$ we have~$F_0=1$.

The image of the centre of the algebra $\YN$ with respect to 
the evaluation homomorphism $\pi$ coincides with $\ZN$. 
Indeed, the images
$\pi(D_1)\lc\pi(D_N)$
of the coefficients of the quantum determinant $D(u)$ generate 
the centre $\ZN$ of the algebra $\UN$; see [C,HU]. 
So we can choose a central element of $\YN$ $$
B_{d_1,d_2,\ldots}^{\ts(s)}\in\pi\1\bigl(F_s\bigr). $$
Moreover, we can assume that for any $s<d_1+2\ts d_2+\ldots$ 
this element is zero. If $d_1=d_2=\ldots=0$ and $s=0$ we can 
assume that this element is $1$.
Now set
$$
\Cal D(z)=\sum_{s\ge0}\ z^{-s}\!\sum_{d_1,\ts d_2,\ts\ldots} 
\
B_{d_1,d_2,\ldots}^{\ts(s)}
\otimes
A_1^{d_1}A_2^{d_2}\ldots
$$
where the inner sum is taken over all sequences 
$d_1,d_2,\ldots$ with $d_1+2\ts d_2+\ldots\le s$.

Consider $f_{\la}^{(s)}$
with fixed $d_1,d_2,\ldots$ as a polynomial in $\la_1\lc\la_N$. 
To prove that
$\pi_\la\ot\id\,\bigr(\Cal D(z)\bigl)=\id\ot D_\la(z)$ for any diagram 
$\la$ it remains to show that this polynomial vanishes when the diagram 
$\la$ satisfies the condition $n<d$; see \(3.9999).

Observe that the coefficient $f_{\la}^{(s)}$ of the series \(3.999) 
has the form of the sum
$$
\sum_{0\le k\le d}
\
(n-k)\ts\ldots\ts(n-d\ts)\cdot g_{\la}^{(k)} $$
where $g_{\la}^{(k)}$ depends on $\la$ as a certain symmetric 
polynomial in $c_1\lc c_n$ and belongs to the ideal generated by 
the $k$-th {\it elementary} symmetric polynomial
$$
\sum_{1\le i_1<\ldots<i_k\le n}
\ts
c_{i_1}\ldots c_{i_k}
\ts.
$$
Let us express the latter sum as a symmetric polynomial in 
$\la_1-1\lc\la_N-N$. We will denote this polynomial by 
$b_k(\la_1\lc\la_N)$. 
It suffices to demonstrate that $b_k(\la_1\lc\la_N)=0$ 
when the diagram $\la$ satisfies the condition $n<k$. 
We will do that by the induction on the difference $k-n=1,2,\ldots\ts$. 

If $\la$ is a diagram with $n\ge1$
boxes then $b_n(\la_1\lc\la_N)=c_1\ldots c_n=0$. 
Indeed, 
then there is at least one box on the main diagonal of the diagram $\la$. 
The content of this box is zero. Further, suppose that $\la$ 
is a diagram such that $\la_1>\ldots>\la_N$. Then due to
$$
\sum_{0\le k\le n}
b_k(\la_1\lc\la_N)
\,
u^{n-k}
\,\ts=
\prod_{1\le k\le n}
(u+c_k)
$$
we have for any $i=1\lc N$ and $k=1\lc n$ the relation 
$$
b_k(\la_1\lc\la_N)=b_k(\la_1\lc\la_i+1\lc\la_N) -(\la_i-i+1)\cdot 
b_{k-1}(\la_1\lc\la_N). 
$$
With the fixed index $k$ this relation must be valid for any 
$\la_1\lc\la_N\in\CC$. In particular, 
if $\la$ is the empty diagram then $b_1(\la_1\lc\la_N)=0$.
The latter fact provides the base for our induction. 
Furthermore, we have already established that $b_k(\la_1\lc\la_N)=0$ 
when $n=k\ge1$.
Now for each $k\ge 2$ the above relation along with 
the inductive assumption implies that $b_k(\la_1\lc\la_N)=0$ 
for any Young diagram $\la$ with less than $k$ boxes
\enddemos

\newpage

Now let $m$ be an arbitrary positive integer and $\mu$ be any 
Young diagram with $m$ boxes.
Equip the set of all pairs $(k,l)$ 
where $1\le k\le n$ and $1\le l\le m$ with the lexicographical ordering.
Let $d_l$ be the content
of the box with $l$ in the column tableaux $\Muc$ of shape $\mu$. 
Consider the rational function of $z$
$$
R_{\la\mu}(z)=
\bigl(F_\la\ot\id\bigr)
\cdot
\prod_{(k,l)}^\rightarrow\ts
R_{\ts k,l+n\ts}(c_k\ts,d_l+z)
\cdot
\bigl(\id\ot F_\mu\bigr)
$$
valued in $\EndCN^{\ot(n+m)}$
where the factors corresponding to the pairs $(k,l)$ 
are arranged with respect to the above ordering. 
By Theorem 2.2 and by the relations \(2.1),\(2.2) any value of 
$R_{\la\mu}(z)$ is divisible on $F_\la\ot\id$ on the 
right and on $\id\ot F_\mu$ on the left. 
Thus we can regard $R_{\la\mu}(z)$ as a function valued in 
$\End(V_\la)\ot\End(V_\mu)$. 
We have the following corollary to Proposition 3.1 and Theorem 3.6. 

\proclaim{Corollary 3.7}
We have the equality 
$\,\pi_\la\ot\pi_\mu\,\bigr(\R(z)\bigl)=R_{\la\mu}(z)$. 
\endproclaim

\nt
In the next section we will be considering behavior of the 
function $R_{\la\mu}(z)$ at $z=0$.

\kern20pt
\centerline{\S\ts\bf4.\ The estimation theorem}\section{\,} \kern-20pt

\nt
Let any two Young diagrams $\la$ and $\mu$ consisting respectively of 
$n$ and
$m$
boxes be fixed.
As before, we will denote by $c_k$ and $d_k$ the contents of the 
boxes with $k$
in the column tableaux $\Lac$ and $\Muc$ respectively. 
We will employ the elements $\Phi_\la\in\CSn$ and $\Phi_\mu\in\CSm$ 
determined by~\(2.0).
%We will keep using Lemma 2.1 and Lemma 2.11. 
Let us equip the set of all pairs $(i,j)$ 
where $1\le i\le n$ and $1\le j\le m$ with the 
lexicographical ordering:
$$
(1,1)\prec\ldots\prec(1,m)
\prec\ldots\ \ldots\ \ldots\prec
(n,1)\prec\ldots\prec(n,m).
\Tag{4.1}
$$
We will regard $S_n$ as a subgroup in $S_{n+m}$ with respect to 
the standard embedding. The superscript $^\vee$ will 
denote the embedding of the group $S_m$ into $S_{n+m}$ determined by 
assigning to the transposition of $j$ and $k$ that of $j+n$ and $k+n$. 
Let $z$ be a complex parameter. Consider the rational function of $z$ 
valued in the group ring $\CSlm$
$$
\Phi_{\la\mu}(z)=\Phi_\la\cdot
\prod_{(i,j)}^\rightarrow\ts
\phi_{\ts i,j+n\ts}(c_i\ts,d_j+z)
\cdot\Phi_\mu^\vee
\Tag{4.2}
$$
where the factors corresponding to the pairs $(i,j)$ are arranged 
with respect to the ordering \(4.1). Let $r$ be the number of the 
boxes on the main diagonal of the diagram $\la$. This number 
is called the {\it rank} of the diagram $\la$.

\proclaim{Theorem 4.1}
The order of the pole of $\Phi_{\la\mu}(z)$ at $z=0$ does not exceed 
the rank $r$ of $\la$.
This order does not exceed $r-1$ if the diagram $\la$ 
is not contained in $\mu$.
\endproclaim

\demo{Proof}
When estimating order of the pole of $\Phi_{\la\mu}(z)$ at $z=0$ 
we employ arguments similar to those already used in the proof of 
Proposition 2.12 but more elaborated. 

The factor $\phi_{\ts i,n+j\ts}(c_i,d_j+z)$ in \(4.2) 
has a pole at $z=0$ if and only if the numbers $i$ and $j$ 
stand on the same diagonals of the tableaux $\Lac$ and $\Muc$ 
respectively. That pole is simple. 
We will then call the pair $(i,j)$ {\it singular}. 
Let us confine every such pair to a certain segment of the 
sequence \(4.1) as follows. 

\smallskip

\vbox{
\nt i)
If $j=1$ then the pair $(i,j)$ itself makes a segment. 

\nt ii)
If $j\neq1$ and $j$ is not in the first row of $\Muc$ then the two pairs 
$(i,j-1)\ts,(i,j)$ form a segment. Note that $j-1$ stands straight above 
$j$ in the tableau $\Muc$.

\nt iii)
If $j\neq1$ appears in the first row of $\Muc$ then the 
corresponding segment is
$$
(i,k)\ts,(i,k+1)\ts\lc(i,j-1)\ts,(i,j)
\Tag{4.25}
$$
where the number $k$ is next to the left of the number $j$ 
in the first row of $\Muc$.}
\smallskip

\nt
Within each of these segments $(i,j)$ is the only one singular pair, 
and all segments
do not intersect. Note that the number of all segments of type (i) is 
exactly $r$.

Consider the segment (ii) where the singular pair $(i,j)$ is fixed. 
Here we have $c_i=d_j=d_{j-1}-1$. By Corollary 2.8 the element 
$\Phi_\mu$ of the group ring $\CSm$ is divisible on the left by 
$\phi_{\ts j-1,j}(d_{j-1}\ts,d_j)$.
Therefore due to the relations \(2.1) and \(2.2) in $\CSlm$ the 
rational function of $z$
$$
\prod_{(i,j)\prec(p,q)}^\rightarrow\ts
\phi_{\ts p,q+n\ts}(c_p\ts,d_q+z)
\cdot\Phi_\mu^\vee
\Tag{4.3}
$$
is divisible on the left by $\phi_{\ts j-1+n,j+n}(d_{j-1}\ts,d_j)$. 
But by Lemma 2.11 the product
$$
\phi_{\ts i,j-1+n}(c_i\ts,d_{j-1}+z)\ts\ts 
\phi_{\ts i,j+n}(c_i\ts,d_j+z)\ts\ts
\phi_{\ts j-1+n,j+n}(d_{j-1}\ts,d_j)
\nopagebreak
$$
is regular at the point $z=0$. 
Thus when estimating from above the order of the pole of \(4.2) at 
$z=0$, the factors of type (ii) do not count. 

Now consider the segment (iii) where the singular pair $(i,j)$ 
is again fixed. Here we have $c_i=d_j=d_{k}+1$. 
By Corollary 2.5 the element $\Phi_\mu$ of $\CSm$ is 
divisible on the left by
$$
\prod_{k<q<j}^\leftarrow\ts
\phi_{\ts q\ts j\ts}(d_q\ts,d_j)
\cdot
\phi_{\ts k\ts j\ts}(d_k\ts,d_j).
$$
Therefore due to the relations \(2.1) and \(2.2) 
the rational function determined by \(4.3) is divisible on the left by 
$$
\prod_{k<q<j}^\leftarrow\ts
\phi_{\ts q+n,j+n\ts}(d_q\ts,d_j)
\cdot
\phi_{\ts k+n,j+n\ts}(d_k\ts,d_j)
\Tag{4.4}
$$
Consider the product of the factors in \(4.2) corresponding to the 
pairs from the segment (iii). Multiplying this product on the right 
by \(4.4) we get 
$$
\gather
\phi_{\ts i,k+n\ts}(c_i\ts,d_k+z)
\cdot\!\!
\prod_{k<q<j}^\rightarrow\!
\phi_{\ts i,q+n\ts}(c_i\ts,d_q+z)
\cdot
\phi_{\ts i,j+n\ts}(c_i\ts,d_j+z)\ts\times \\
\prod_{k<q<j}^\leftarrow\!
\phi_{\ts q+n,j+n\ts}(d_q\ts,d_j)
\cdot
\phi_{\ts k+n,j+n\ts}(d_k\ts,d_j)\,\,\,= \prod_{k<q<j}^\leftarrow\!
\phi_{\ts q+n,j+n\ts}(d_q\ts,d_j)\ts\times \\
\phi_{\ts i,k+n\ts}(c_i\ts,d_k+z)\ts\ts
\phi_{\ts i,j+n\ts}(c_i\ts,d_j+z)\ts\ts
\phi_{\ts k+n,j+n\ts}(d_k\ts,d_j)\ts
\cdot\!\!
\prod_{k<q<j}^\rightarrow\!
\phi_{\ts i,q+n\ts}(c_i\ts,d_q+z)
\endgather
$$
where we again used the relations \(2.1)\ts,\(2.2). 
The rational function at the right hand side of the latter equality 
is regular at $z=0$ by Lemma 2.11. 
Thus when estimating the order of the pole of \(4.2) at $z=0$, 
the factors of the type (iii) do not count either. This completes
the proof of the first statement in Theorem 4.1. \Par
\vbox{

The proof of the second statement is similar but more involved. 
As before, denote
\line{by $\laps\hskip-0.4pt$ and $\mups$ the lengths of the 
$\hskip-0.5pt s\hskip-0.4pt$-th columns of the diagrams $\la$ and 
$\mu$ respectively.}
Assume that the diagram $\la$ is not contained in $\mu$. 
We shall demonstrate that when estimating the order of the pole of 
\(4.2) at $z=0$, we may exclude from counting the factor corresponding 
to the singular pair $(1,1)$ as well as those corresponding to 
the singular pairs $(i,j)$~with~$i\neq1$. }
Fix the minimal
number $s$ such that $\laps>\mups\ts$. Let $i_1\lc i_s$ and $j_1\lc j_s$ 
be the first $s$ numbers in the first rows of the tableaux $\Lac$ and
$\Muc$ respectively. Of course, here $i_1=j_1=1$. Denote $\mups-1=t$ 
so that the numbers in the $s\hskip-0.3016pt$-th column of $\Muc$ are 
$j_s,j_s+1\lc j_s+t$. 
The singular pairs
$$
(i_1\ts,j_1)\lc(i_s\ts,j_s)
\quad\text{and}\quad
(i_s+1\ts,j_s+1)\lc(i_s+t\ts,j_s+t)
$$
will be called {\it special\,.} 
We will now proceed in several steps. 

I.
Let us reorder the pairs
$(i,j)$ in the product \(4.2) as follows\,; 
this reordering will not affect the value of the product due 
to the relations \(2.2). Namely, for each $i=i_s\ts,i_s+1\lc i_s+t$ 
change the sequence of pairs 
$$
(i,j)\ts,(i,j+1)\lc(i,m)\ts,(i+1,1)\lc(i+1,j-1)\ts,(i+1,j) 
$$
where $j=j_s\ts,j_s+1\lc j_s+t$ respectively, to the sequence 
$$
(i+1,1)\lc(i+1,j-1)\ts,(i,j)\ts,(i+1,j)\ts,(i,j+1)\lc(i,m). 
$$
The latter reordering preserves all the segments (i-iii) except 
for the segment with $(i_s,j_s)$ and the $t$ segments of type (ii) 
$$
\align
&(i_s+q,j_s+q-1)\ts,(i_s+q,j_s+q)\,;
\qquad
q=1\lc t\,.
\\
\intertext{Instead of these $t+1$ 
segments we will introduce the segments in our new ordering
}
&(i_s+q,j_s+q)\ts,(i_s+q+1,j_s+q)\,;
\qquad
q=0\ts,1\lc t\,.
\Tag{4.5}
\endalign
$$
Observe that here
$$
d_{j_s+q}=c_{i_s+q}=c_{i_s+q+1}+1.
$$
Moreover, by Corollary 2.8
the element $\Phi_\la$ of $\CSn$ is divisible on the right by 
$$
\phi_{\ts i_s+q,i_s+q+1\ts}(c_{i_s+q\ts},c_{i_s+q+1}). 
$$
As well as in the proof of the first statement of Theorem 4.1, 
but by making now use of
Lemma 2.1 we can now exclude from our counting the factors in 
\(4.2) corresponding to each of the special singular pairs 
$$
(i_s+q\ts,j_s+q)\ts;
\qquad
q=0\ts,1\lc t\,.
\nopagebreak
$$
{If $\mups=0$ this step is empty. If $s=1$ then the proof of 
Theorem 4.1 is completed.}

II. Suppose that $s>1$. Note that by our assumption $\lapsmin\le\mupsmin$. 
Denote
$$
i_s=i,
\quad
j_s=j,
\quad
i_{s-1}=h,
\quad
j_{s-1}=k
\quad\text{and}\quad
j_{s-1}+\lapsmin-1=l.
$$
On the left of the next picture we show the $(s-1)$-th and $s$-th 
columns of the tableau $\Lac$ while on the right we have the same 
columns of the tableau $\Muc$.

\smallskip\smallskip\smallskip

\vsk>
\vbox{
$$
{\bx}
{\bx}
\phantom{\bx}
\phantom{\bx}
\phantom{\bx}
\phantom{\bx}
{\bx}
{\bx}
$$
\vglue-17.8pt
$$
{\bx}
{\bx}
\phantom{\bx}
\phantom{\bx}
\phantom{\bx}
\phantom{\bx}
{\bx}
{\bx}
$$
\vglue-17.8pt
$$
{\bx}
{\bx}
\phantom{\bx}
\phantom{\bx}
\phantom{\bx}
\phantom{\bx}
{\bx}
{\bx}
$$
\vglue-17.8pt
$$
{\bx}
{\bx}
\phantom{\bx}
\phantom{\bx}
\phantom{\bx}
\phantom{\bx}
{\bx}
\phantom{\bx}
$$
\vglue-17.8pt
$$
{\bx}
\phantom{\bx}
\phantom{\bx}
\phantom{\bx}
\phantom{\bx}
\phantom{\bx}
{\bx}
\phantom{\bx}
$$
\vglue-17.8pt
$$
{\bx}
\phantom{\bx}
\phantom{\bx}
\phantom{\bx}
\phantom{\bx}
\phantom{\bx}
{\bx}
\phantom{\bx}
$$
\vglue-17.8pt
$$
\phantom{\bx}
\phantom{\bx}
\phantom{\bx}
\phantom{\bx}
\phantom{\bx}
\phantom{\bx}
{\bx}
\phantom{\bx}
$$
\vglue-17.8pt
$$
\phantom{\bx}
\phantom{\bx}
\phantom{\bx}
\phantom{\bx}
\phantom{\bx}
\phantom{\bx}
{\bx}
\phantom{\bx}
$$
\vglue-124pt
$$
h
\kern10pt
i
\kern67pt
k
\kern9pt
j
$$
\vglue-18pt
$$
\kern1pt
\cdot
\kern9pt
\cdot
\kern66pt
\cdot
\kern12pt
\cdot
$$
\vglue-18pt
$$
\kern1pt
\cdot
\kern9pt
\cdot
\kern66pt
\cdot
\kern12pt
\cdot
$$
\vglue-18pt
$$
\kern1pt
\cdot
\kern9pt
\cdot
\kern66pt
\cdot
\kern12pt
\phantom{\cdot}
$$
\vglue-18pt
$$
\kern3pt
\cdot
\kern10pt
\phantom{\cdot}
\kern67pt
\cdot
\kern12pt
\phantom{\cdot}
$$
\vglue-17.5pt
$$
\kern2pt
\cdot
\kern10pt
\phantom{\cdot}
\kern68pt
l
\kern12pt
\phantom{\cdot}
$$
\vglue-18pt
$$
\kern71pt\cdot
$$
\vglue-17.5pt
$$
\kern71pt\cdot
\nopagebreak
$$
}
\vsk>

\smallskip

\nt
Observe
that in this notation the segment of $(i,j)$ in \(4.1) was exactly \(4.25). 

Let us perform further reordering of the factors in \(4.2). Consider the
subsequence in \(4.1) consisting of all the pairs $(p,q)$ such that 
$$
(h,k)\preccurlyeq(p,q)\preccurlyeq(i,l)\ts; 
$$
this subsequence has been not affected by reordering at the previous step.
Let us now consecutively take the pairs $(p,q)$ with $h\le p<i$ and 
$n<q\le m$ to the left of that subsequence. Then take the pairs 
$(p,q)$ with $h<p\le i$ and $1\le q<k$ to the right. 
Between them will remain the sequence
$$
(h,k)\ts\lc(h,l)\ts,(h+1,k)\ts\lc(h+1,l)\ts\lc\ldots\lc(i,k)\ts\lc(i,l)\ts.
$$
This reordering does not alter the value of \(4.2) again due to the 
relations \(2.2). It does not break any of the new segments \(4.5). 
It also preserves all the segments \text{(i-iii)}
except for those containing the singular pairs $(i,j)$ and $(h,k)$. 
The factor in \(4.2) corresponding to $(i,j)$ has been excluded from our 
counting at the previous step. By Corollary 2.5
the element $\Phi_\la$ of $\CSn$ is divisible on the right by 
$$
\phi_{\ts hi\ts}(c_h\ts,c_i)\ts\ts\cdot\hskip-2pt 
\prod_{h<p<i}^\rightarrow\ts
\phi_{\ts pi\ts}(c_p\ts,c_i).
$$
Note that here $d_k=c_h=c_i-1$.
Therefore by using Lemma 2.1 along with the relations \(2.1)\ts,\(2.2) 
we can exclude from our counting the factor in \(4.2) corresponding to 
the pair $(h,k)$ as well. This ends the second step. 

Let us now repeat the last step for $s-1\lc 2$ instead of $s$. Then all 
the factors in \(4.2) corresponding to special singular pairs will be 
excluded from our counting. In particular, the factor corresponding to 
$(1,1)$ will be excluded.
The segments \text{(i-iii)} with non-special singular pairs will be 
preserved by every reordering. So we will exclude from counting the 
factors in \(4.2) corresponding to all non-special singular pairs 
$(i,j)$ with $j\neq 1$.
This can be demonstrated by arguments used in the proof of the 
first statement of Theorem 4.1. The proof of the second statement is 
now also completed
\enddemos

\kern20pt
\centerline{\S\ts\bf5.\ Higher Capelli identities}\section{\,} 
\kern-20pt

\nt
Let any two positive integers $N$ and $M$ be fixed. In this section 
we will consider invariants of the action of Lie algebra $\glMN$ in 
the
space $\PDMN$ of differential operators on $\CMN$ with polynomial 
coefficients. Let the indices $i,j$ and $a,b$ run through $1\lc N$ and 
$1\lc M$ respectively.
Let $x_{ia}$ be the standard coordinates on the vector space $\CMN$. 
We will define the actions of $\glN$ and $\glM$ in the space $\PMN$ 
of polynomial functions on $\CMN$ by 
%the assignments
$$
E_{ij}\mapsto\sum_b\,x_{ib}\,\d_{jb}
\qquad\text{and}\qquad
E_{ab}\mapsto\sum_j\,x_{ja}\ts\d_{jb}
\Tag{5.1}
$$
respectively. Here $\d_{jb}$ is the partial derivation with respect to 
the coordinate $x_{jb}$.
We have the irreducible decomposition [HU] of the $\glMN$-module $$
\PMN\,=\,
\underset\la\to\oplus\,\ V_\la\ot W_\la
\Tag{5.0}
$$
where $\la$ runs through the set of all Young diagrams with $0,1,2,\ldots$ 
boxes and not more than $M,N$ rows. Here $W_\la$ is the irreducible 
$\glM$-module corresponding to the diagram $\la$. 
Thus the action of $\glN\times\glM$
in $\PMN$ is multiplicity-free.

The space $\DMN$ of differential operators on $\CMN$ with constant 
coefficients can be regarded as the $\glMN$-module dual to $\PMN$. 
Hence the space $\Cal I$ of all the $\glMN$-invariants in
$$
\PDMN=\PMN\ot\DMN
\Tag{5.2}
$$
is a direct sum of one-dimensional subspaces corresponding to the 
diagrams $\la$ with not more than $M,N$ rows.
Denote by $\Cal I_\la$ the one-dimensional
subspace corresponding to the partition $\la$. Let us give an explicit 
formula for a non-zero vector from $\Cal I_\la$. Suppose that the 
diagram $\la$ consists of $n$ boxes. Let
$\chi_{\ts\la}$ be the character of the irreducible $S_n$-module $U_\la$. 
Consider the element of $\PDMN$ defined by \(1.1). We will denote
this element by $c_\la$. 
%$$
%c_{\ts\la}\,=
%\sum_{\si\ts\in\,S_n}
%\ \,
%\sum_{i_1\lc i_n}
%\,
%\sum_{a_1\lc a_n}
%\
%{\chi_\la\ts(\si)}
%\cdot
%x_{i_1a_1}\ldots x_{i_na_n}
%\!\!\cdot
%\d_{i_{\si(1)}a_1}\ldots\d_{i_{\si(n)}a_n} 
%$$
%where the indices $i_1\lc i_n$ and $a_1\lc a_n$ run 
%through $1\lc N$ and $1\lc M$.
%respectively.
%This element is invariant under the action of $\glMN$. 

\proclaim{Proposition 5.1}
We have $c_\la\in\Cal I_\la$.
\endproclaim

\demo{Proof}
Consider the $\glMN$-submodule in $\PMN$ formed by polynomials of the 
degree $n$. This submodule can be identified with the subspace of 
invariants in $(\CC^N)^{\ot n}\ot(\CC^M)^{\ot n}$ with respect to 
simultaneos permutations of the first $n$ and the last
$n$ tensor factors. The permutation $\si$ of the first $n$ tensor 
factors
corresponds to the map
$$
x_{i_1a_1}\ldots x_{i_na_n}\,\mapsto\,
x_{i_{\si(1)}a_1}\ldots x_{i_{\si(n)}a_n} 
\nopagebreak
$$
in $\PMN$. We have the irreducible decomposition of the 
$S_n\times\glN$-module
$$
(\CC^N)^{\ot n}\,=\,
\underset\la\to\oplus\,\ U_\la\ot V_\la
$$
where $\la$ runs through the set of all diagrams with $n$ boxes 
and not more than $N$ rows. So the image of the map 
$$
\eta\,:\ \,
x_{i_1a_1}\ldots x_{i_na_n}\,\mapsto\,
\sum_{\si\ts\in\,S_n}
\frac{\chi_\la\ts(\si)\!}{n!}
\cdot
x_{i_{\si^{-1}(1)}a_1}\ldots\,x_{i_{\si^{-1}}(n)a_n} $$
in $\PMN$ is the subspace $V_\la\ot W_\la$ if the diagram $\la$ has no 
more
than $M,N$ boxes.
The element $c_\la$ in \(5.2) is the image with respect to the map 
$\eta\ot\id$ of $$
\sum_{i_1\lc i_n}
\,
\sum_{a_1\lc a_n}
\
x_{i_1a_1}\ldots x_{i_na_n}
\!\!\cdot
\d_{i_1a_1}\ldots\d_{i_na_n}.
$$ The latter element of $\PDMN$ is manifestly $\glMN$-invariant 
\enddemos

\nt
Extend the action \(5.1) of the Lie algebra $\glN$ in $\PMN$ to 
the universal enveloping algebra $\UN$.
The subspace of $\glMN$-invariants
in $\PDMN$ coincides with the image of the centre $\ZN$ of the algebra 
$\UN$;
see [HU]. We will point out an element $e_\la\in\ZN$ whose image is 
$c_\la$.
Let $\psi_\la:\YN\to\CC$ be the trace of the representation $\pi_\la$.
Consider~the~element
$$
\psi_\la\ot\pi\bigl(\R(z)\bigr)\in\UN\,[[z\1]] $$
where $\pi$ is the evaluation homomorphism \(3.5). Then due to Theorem 3.6 
the product \(1.2) is a polynomial in $z$. 
Denote by $e_\la(z)$ this polynomial.
The following general observation goes back to [D2]. 

\proclaim{Proposition 5.2}
The coefficients of the polynomial $e_\la(z)$ belong to $\ZN$. 
\endproclaim

\demo{Proof}
\!Consider the universal enveloping algebra $\UN$ as a subalgebra in 
$\!\YN$.
Every element of this subalgebra is invariant with respect to the 
automorphism
$\tau_z$ for any $z\in\CC$.
Furthermore, restriction of the comultiplication \(3.3) to $\UN$ is 
cocommutative:
$$
\De(Y)=Y\!\ot1+1\ot Y,
\qquad
Y\in\glN\ts.
$$
Therefore by the definition of the series $\R(z)$ for any $Y\in\glN$ 
the commutator
$$
\bigl[\ts\R(z)\,,Y\!\ot1+1\ot Y\ts\bigr]=0. $$
By applying the map $\psi_\la\ot\pi$ to the latter equality we obtain 
that
$$
\bigl[\ts\psi_\la\ot\pi\bigl(\R(z)\bigr)\,,Y\ts\bigr]=0\quad\square $$
\enddemo

\nt
The next theorem is the main result of this section. Denote by $e_\la$
the value $e_\la(0)\hskip-1pt$.

\proclaim{Theorem 5.3}
The image in $\PDMN$ of the element $e_\la\in\ZN$~is~$c_\la.$ 
\endproclaim

\nt
We will present the main steps of the proof as separate propositions. 
First let us give an explicit formula for the polynomial in $z$ 
valued in $\End(V_\la)\ot\UN$
$$ 
E_\la(z)=\pi_\la\ot\pi\bigl(\R(z)\bigr)\cdot(z-c_1)\ldots(z-c_n) 
\Tag{5.99}
$$
As well as before, here we regard $\End(V_\la)$ as a subalgebra in 
$\EndCN^{\ot n}$.

\nt
Introduce the element
$$
E=\ts\sum_{i,j}\ts E_{ji}\ot E_{ij}\in\EndCN\ot\UN\ts 
$$
where the indices $i,j$ run through $1\lc N$. 
%Put $E(u)=u+E$.
For $s=1\lc n$ we will write
$$
E_s=\iota_s\ot\id\ts\bigl(E\bigr)\in\EndCN^{\ot n}\ot\UN. 
$$
Then due to Theorem 3.6 and the definition \(3.5) 
$$
\align
E_\la(z)
&=\id\ot\pi\bigl(T_\la(z)\bigr)\cdot(z-c_1)\ldots(z-c_n) 
\Tag{5.55}
\\
&=F_\la\ot1\cdot(E_1+z-c_1)\ldots(E_n+z-c_n). \endalign
$$

Let $\tr$ denote the standard matrix trace on $\EndCN^{\ot n}$. 
By definition we have
$$
e_\la(z)=\tr\ot\id\ts\bigl(E_\la(z)\bigr). \Tag{5.56}
$$
Denote by $D$ the element of the algebra $\EndCN^{\ot n}\ot\PDMN$ 
$$
\sum_{i_1\lc j_n}\,
\sum_{j_1\lc j_n}
\sum_{a_1\lc a_n}
\
E_{j_1i_1}\!\ot\ldots\ot E_{j_ni_n}
\ot
x_{i_1a_1}\ldots x_{i_na_n}
\d_{j_1a_1}\ldots\d_{j_na_n}.
$$ 
Consider the product $C_\la=F_\la\ot1\cdot D.$ 
We will employ the following observation. 

\proclaim{Proposition 5.4}
We have the equality
$c_\la=\tr\ot\id\ts(C_\la).$
\endproclaim

\demo{Proof}
For any $F$ from the image of the symmetric group $S_n$ in 
$\EndCN^{\ot n}$ the elements $F\ot1$ and $D$ of the algebra 
$\EndCN^{\ot n}\ot\PDMN$ commute. Thus Proposition 5.4 follows from the
definition of $c_\la$ and the identity in $\CSn$ $$
\frac1{\ts n!}\,
\sum_{\si\in S_n}\,
\si\,\ts\Phi_\la\,\si\1=
\frac1{\ts\,\dim U_\la}
\cdot\sum_{\si\in S_n}\,
{\chi_\la\ts(\si)}
\,\ts
\si\quad\square
$$
\enddemo

\nt
By Proposition 5.4 and the equality \(5.56) the next proposition implies 
Theorem~5.3.

\proclaim{Proposition 5.5}
The image of $E_\la(0)$ in $\EndCN^{\ot n}\ot\PDMN$ is $C_\la$. 
\endproclaim

\demo{Proof}
Consider the standard gradation on the vector space $\PDMN$ by the 
degree of a differential operator. Extend this gradation to the vector 
space $\EndCN^{\ot n}\ot\PDMN$. Denote the image of $E_\la(0)$
in that space by~$C_\la^{\,\prime}$. Due to \(5.55) the leading term of 
$C_\la^{\,\prime}$ coincides with the homogeneous element $C_\la$. 
Thus it remains to show that the element $C_\la^{\,\prime}$ is 
homogeneous as well. Due to \(5.0) it suffices to check that for any 
Young diagram $\mu$ with less than $n$ boxes $$
\id\ot\pi_\mu\bigl(E_\la(0)\bigr)=0.
\Tag{5.98}
$$
In this case the diagram $\la$ is not contained in $\mu$. 
The number of zeroes among~the contents $c_1\lc c_n$ of the 
boxes of the diagram $\la$ is exactly its rank $r$. So by the definition 
\(5.99)
the equality \(5.98) follows from Corollary 3.7 and Theorem~4.1 
\enddemos

\kern15pt
\centerline{\bf Acknowledgements}\section{\,} \kern-20pt

\nt
The author has been supported by the EPSRC Advanced Research Fellowship.
He was also supported by the Erwin Schr\"odinger International Institute 
in Vienna.

\newpage
\vbox{
\centerline{\bf References}\section{\,}
\kern-20pt

\itemitem{[C]}
{A. Capelli},
{\it Sur les op\'erations dans la th\'eorie des formes alg\'ebriques},
{Math. Ann.}
{\bf 37}
(1890),
1--37.

\itemitem{[C1]}
{I. Cherednik},
{\it On special bases of irreducible finite-dimensional representations of 
the degenerate affine Hecke algebra},
{Funct. Anal. Appl.}
{\bf 20}
(1986),
87--89.

\itemitem{[C2]}
{I. Cherednik},
{\it A new interpretation of Gelfand-Zetlin bases}, {Duke Math. J.}
{\bf 54}
(1987),
563--577.

\itemitem{[D1]}
{V. Drinfeld},
{\it Hopf algebras and the quantum Yang--Baxter equation}, {\text{Soviet} 
Math.\,Dokl.}
{\bf 32}
(1985),
254--258.

\itemitem{[D2]}
{V. Drinfeld},
{\it On almost cocommutative Hopf algebras}, {Leningrad Math. J.}
{\bf 1}
(1990),
321--342.

\itemitem{[GM]}
{A. Garsia and T. Maclarnan},
{\it Relations between Young's natural and the Kazhdan-Lusztig 
representations of $S_n$},
{Adv. Math.}
{\bf 69}
(1988),
32--92.

\itemitem{[HU]}
{R. Howe and T. Umeda},
{\it The Capelli identity, the double commutant theorem, and 
multiplicity-free actions},
{Math. Ann.}
{\bf 290}
(1991),
569--619.

\itemitem{[J]}
{A. Jucys},
{\it On the Young operators of the symmetric groups}, 
Lietuvos Fizikos Rinkinys
{\bf 6}
(1966),
163--180.

\itemitem{[JKMO]}
M. Jimbo, A. Kuniba, T. Miwa and M. Okado, 
{\it The $A_n^{(1)}$ face models},
Comm. Math. Phys.
{\bf 119}
(1988),
543--565.

%\itemitem{[KR]}
%{A. Kirillov and N. Reshetikhin},
%{\it Yangians, Bethe ansatz and combinatorics}, 
%{Lett.\,Math.\,Phys.}
%{\bf 12}
%(1986),
%199--208.

\itemitem{[KRS]}
P. Kulish, N. Reshetikhin and E. Sklyanin, 
{\it Yang--Baxter equation and the representation theory I}, 
Lett. Math. Phys.
{\bf 5}
(1981),
393--403.

\itemitem{[MNO]}
{A. Molev, M. Nazarov and G. Olshanski}, 
{\it Yangians and classical Lie algebras}, Russian Math. Surveys
{\bf 51}:2
(1996),
27--104.

\itemitem{[N1]}
{M. Nazarov},
{\it Quantum Berezinian and the classical Capelli identity}, 
{Lett. Math. Phys.}
{\bf 21}
(1991),
123--131.

\itemitem{[N2]}
{M. Nazarov},
{\it Young's symmetrizers for projective representations of the 
symmetric group}, preprint;
%{Preprint RIMS},
%{\bf 900},
%Kyoto, 1992;
to appear in {Adv. Math.}

\itemitem{[O]}
{G. Olshanski},
{\it Representations of infinite-dimensional classical groups,
limits of enveloping algebras, and Yangians,} in
{\lq\lq\ts Topics in Representation Theory\ts\rq\rq}, 
%(A. A. Kirillov, Ed.)
{Adv. Soviet Math.}
{\bf 2},
Amer. Math. Soc.,
Providence,
1991,
pp. 1--66.

\itemitem{[O1]}
{A. Okounkov},
{\it Quantum immanants and higher Capelli identities}, 
Transformation Groups,
{\bf 1}
(1996),
99--126.

\itemitem{[O2]}
{A.\,Okounkov},
{Young basis, Wick formula, and higher Capelli identities}, 
pre\-print; to appear in Int. Math. Research Notes.

\itemitem{[OO]}
{A. Okounkov and G. Olshanski},
{\it Shifted Schur functions},
preprint. 
%to appear in S.-Petersburg Math. J.

\itemitem{[RS]}
N. Reshetikhin and M. Semenov-Tian-Shansky, 
{\it Central extensions of quantum current groups}, Lett. Math. Phys.
{\bf 19}
(1990),
133--142.

\itemitem{[S]}
{S. Sahi},
{\it The spectrum of certain invariant differential operators associated 
to a Hermitian symmetric space,}
in
{\lq\lq\ts Lie Theory and Geometry\ts\rq\rq}, 
%(J.-L. Brylinski, R. Brylinski, V. Guillemin, V. Kac, Eds.) 
{Progress in Math.}
{\bf 123},
Birkh\"auser,
Boston,
1994,
pp. 569--576.

%\itemitem{[W]}
%H. Weyl,
%{\it Classical groups, their invariants and representations}, 
%Pri\-nce\-ton University Press,
%Princeton,
%1946.

\itemitem{[Y1]}
{A. Young,}
{\it On quantitative substitutional analysis I} and {\it II}, 
{Proc. London Math. Soc.}
{\bf 33}
(1901),
97--146
{and}
{\bf 34}
(1902),
361--397.

\itemitem{[Y2]}
{A.\, Young,}
{\it On quantitative substitutional analysis VI,} 
{Proc.\,London Math. Soc.}
{\bf 31}
(1931),
253--289.

\bigskip

\nt
\line{\smc
Mathematics Department, University of Wales, Swansea SA2 8PP, UK\hfill}} 

\bye